\documentclass{amsart}
\usepackage{color}

\usepackage{amsmath} 
\usepackage{amssymb}
\usepackage{mathrsfs}

\newtheorem{theorem}{Theorem}[section] 
\newtheorem{claim}[theorem]{Claim}
\newtheorem{mc}[theorem]{Main Claim}

\newtheorem{cd}[theorem]{Claim/Definition}

\newtheorem{conclusion}[theorem]{Conclusion}
\newtheorem{observation}[theorem]{Observation}

\theoremstyle{definition}
\newtheorem{definition}[theorem]{Definition}
\newtheorem{convention}[theorem]{Convention}

\newtheorem{conjecture}[theorem]{Conjecture}
\newtheorem{fact}[theorem]{Fact}

\newtheorem{discussion}[theorem]{Discussion}

\theoremstyle{remark}
\newtheorem{remark}[theorem]{Remark}
\newtheorem{question}[theorem]{Question}
\newtheorem{notation}[theorem]{Notation}

\newcommand{\Th}{{\rm Th}}
\newcommand{\TV}{{\rm TV}}

\newcommand{\SOP}{{\rm SOP}}

\newcommand{\BPA}{{\rm BPA}}

\newcommand{\NSOP}{{\rm NSOP}}

\newcommand{\PA}{{\rm PA}}

\newcommand{\rp}{{\rm rp}}

\newcommand{\lex}{{\rm lex}}
\newcommand{\sym}{{\rm sym}}
\newcommand{\fil}{{\rm fil}}

\newcommand{\true}{{\rm true}}

\newcommand{\feq}{{\rm feq}}
\newcommand{\ceq}{{\rm ceq}}
\newcommand{\fcp}{{\rm fcp}}
\newcommand{\rt}{{\rm rt}}

\newcommand{\tr}{{\rm tr}}
\newcommand{\uf}{{\rm uf}}

\newcommand{\comp}{{\rm comp}}
\newcommand{\false}{{\rm false}}
\newcommand{\iif}{{\rm if}}
\newcommand{\ord}{{\rm ord}}
\newcommand{\Mod}{{\rm Mod}}

\newcommand{\ba}{{\rm ba}}
\newcommand{\pair}{{\rm pair}}

\newcommand{\RSP}{{\rm RSP}}

\newcommand{\id}{{\rm id}}

\newcommand{\Dom}{{\rm Dom}}

\newcommand{\rest}{{\restriction}}
\newcommand{\dom}{{\rm dom}}

\newcommand{\wilog}{{\rm without loss of generality}}
\newcommand{\Wilog}{{\rm Without loss of enerality}}

\def\emptycs{}
\def\evaluateROMANlist{%
        \ifx\ROMANlist\emptycs\else
        \expandafter\xxevaluate\ROMANlist\xxfertig\evaluateROMANlist\fi}
\def\xxevaluate#1,#2\xxfertig{\expandafter\newcommand\csname#1\endcsname{\mathrm
{#1}}\def\ROMANlist{#2}}
 \newcommand{\ROMANlist}{Prop,pseudoProd,nisayon,nisuy,}
 \evaluateROMANlist

\newcommand{\then}{{\underline{then}}}
\newcommand{\when}{{\underline{when}}}

\newcommand{\Then}{{\underline{Then}}}

\newcommand{\If}{{\underline{if}}}
\newcommand{\Iff}{{\underline{iff}}}

\newcommand{\sn}{{\smallskip\noindent}}

\newcommand{\cB}{{\mathscr B}}

\newcommand{\gB}{{\mathfrak B}}

\newcommand{\cC}{{\mathscr C}}

\newcommand{\cL}{{\mathscr L}}

\newcommand{\cI}{{\mathscr I}}

\newcommand{\bbL}{{\mathbb L}}
\newcommand{\bbN}{{\mathbb N}}

\newcommand{\cP}{{\mathscr P}}
\newcommand{\gp}{{\mathfrak p}}

\newcommand{\bbQ}{{\mathbb Q}}

\newcommand{\cS}{{\mathscr S}}
\newcommand{\cT}{{\mathscr T}}
\newcommand{\gt}{{\mathfrak t}}

\newcommand{\cf}{{\rm cf}}

\newcount\skewfactor
\def\mathunderaccent#1#2 {\let\theaccent#1\skewfactor#2
\mathpalette\putaccentunder}
\def\putaccentunder#1#2{\oalign{$#1#2$\crcr\hidewidth
\vbox to.2ex{\hbox{$#1\skew\skewfactor\theaccent{}$}\vss}\hidewidth}}

\renewcommand{\thesubsection}{\arabic{section}(\Alph{subsection})}

\newenvironment{PROOF}[2][\proofname.]
   {\begin{proof}[#1]\renewcommand{\qedsymbol}{\textsquare$_{\rm #2}$}}
   {\end{proof}}

\usepackage[hidelinks]{hyperref}
\begin{document}
\makeatletter\def\shfiuwefootnote{\gdef\@thefnmark{}\@footnotetext}\makeatother\shfiuwefootnote{Version 2023-06-12\_3. See \url{https://shelah.logic.at/papers/1064/} for possible updates.}

\title {Atomic Saturation of reduced powers \\
Sh1064}
\author {Saharon Shelah}
\address{Einstein Institute of Mathematics\\
Edmond J. Safra Campus, Givat Ram\\
The Hebrew University of Jerusalem\\
Jerusalem, 91904, Israel\\
 and \\
 Department of Mathematics\\
 Hill Center - Busch Campus \\ 
 Rutgers, The State University of New Jersey \\
 110 Frelinghuysen Road \\
 Piscataway, NJ 08854-8019 USA}
\email{shelah@math.huji.ac.il}
\urladdr{http://shelah.logic.at}

\thanks{The research of the author was partially supported by European Research Council, Grant No. 338821. The author thanks Alice Leonhardt for the beautiful typing. First typed: February 26, 2014.\\ 
For the changes after publication, the author would like to thank ISF-BSF for partially supporting this research by grant with Maryanthe Malliaris number NSF 2051825, BSF 3013005235.
In the new version, the author is grateful for the generous funding of typing services donated by a person who wishes to remain anonymous and would like to thank the typists for the careful and beautiful typing.  \\ 
References like \cite[Th0.2=Ly5]{Sh:950} means the label of Theorem 0.2 is y5.  The reader should note that the version in my website is usually more updated than the one in the mathematical archive. Was proof read with proofs for the Journal.}  



 
\date{June 11, 2023}

\begin{abstract}
    Our aim was to try to 
    generalize some theorems about the saturation 
    of ultra-powers to reduced powers.  Naturally, we deal 
    with 
    saturation
    for types consisting of atomic formulas.  We succeed to generalize
    ``the theory of dense linear  order 
    (or $ T $ with the strict order property) 
    is maximal and so is any
    $\pair\,(T,\Delta)$ which is $\SOP_3$", (where $\Delta$ consists of atomic or
    conjunction of atomic formulas).  However,
    the theorem on {`}{`}it is enough to deal with symmetric pre-cuts" 
    (so
    the $\gp = \gt$ theorem) cannot 
    be generalized in this case.  Similarly the
    uniqueness of the  dual cofinality fails in this context.

    After publication, a full proof of 2.10m3as added 
    (and minor proofs done). 
\end{abstract}

\maketitle
\numberwithin{equation}{section}
\setcounter{section}{-1}
\newpage

\centerline{Annotated Content}

\S0 \quad Introduction, 
pg.\pageref{intro}

\S1 \quad Axiomatizing \cite[Ch.VI,2.6]{Sh:c}, 
pg. \pageref{1}

\begin{enumerate}
\item[${{}}$]  [We phrase and prove a theorem which 
axiomatizes \cite[Ch.VI,2.6]{Sh:c}.  
The theorem there says that if $D$ is a regular ultra-filter
  on $I$ and for every model $M$ of the theory of dense linear orders
  (or $ T $ with the strict order property),
the model $M^I/D$ is $\lambda^+$-saturated, then $D$ is $\lambda^+$-good and
$\lambda$-regular.] 
\end{enumerate}

\S2 \quad Applying the axiomatized frame, 
pg.\pageref{2}

\begin{enumerate}
\item[${{}}$]  [The axiomatization in \S1 can be phrased as a set of
sentences, surprisingly moreover Horn ones (first order if $\theta_{\mathbf r} =
  \aleph_0$).  Now in this case we can straightforwardly derive
 \cite[Ch.VI,2.6]{Sh:a}.  
 However we can get more, 
 because the axiomatization being Horn,
 we can now deal also with the ($\lambda^+$,atomic)-saturation of
 reduced powers.  We then deal with infinitary logics and comment on models
 of Bounded Peano Arithmetic. After an inquiring of G. Cherlin the proof of \ref{c31} was given.] 
\end{enumerate}

\S3 \quad Criterion for atomic saturation of reduced powers, 
pg.\pageref{3}

\begin{enumerate}
\item[${{}}$]  [For a complete first order $T$ we characterize when a
  filter $D$ on $I$ is such that $M^I/D$ is
($\lambda$,atomic)-saturated for every model $M$ of $T$.]
\end{enumerate}

\S4 \quad The counterexample, 
pg.\pageref{4}

\begin{enumerate}
\item[${{}}$]  [We prove that for reduced powers, the parallel of $\gt
  \le \gp$ in general fails. Also,
  similarly 
  the uniqueness of the dual cofinality.
  More specifically, for $ \lambda \ge {\aleph_1}$, 
  for some regular filter $ D $ on $ \lambda $, the 
  partial order $ (\mathbb{Q} , <  )\lambda /D $  has
  no symmetric pre-cut of cofinality $ \ge \lambda $ 
  but has such an asymmetric pre-cut.]
\end{enumerate}

%
%
\newpage

\section {Introduction} \label{intro}

\subsection {Background, Questions and Answers} \

We know much on saturation of ultrapowers, see Keisler  \cite{Ke67},
\cite[Ch.IV]{Sh:c} and later mainly works of Malliaris and the author,
e.g. (\cite{Sh:998}, \cite{Sh:1030}).  But we know considerably 
less on reduced powers.  For
transparency, let $T$ denote a first-order complete countable theory with
elimination of quantifiers and $M$ will denote a model of $T$.  For
$D$ a regular filter on $\lambda > \aleph_0$ we may ask:
when  is  $M^\lambda/D$
 $ \lambda^+ $-saturated?   
 For $D$ an ultra-filter, Keisler \cite{Kei64} proves that
this holds for every $T$ iff $D$ is $\lambda^+$-good iff this holds
for $T=$ theory of Boolean algebras, such $T$ is called
$\trianglelefteq_\lambda$-maximal.

By \cite[Ch.VI,2.6]{Sh:a} the maximality holds for $T=$ theory of dense
linear orders or just any $T$ with the strict order property 
and by \cite{Sh:500}, any $T$ with the 3-strong order property,
$\SOP_3$ is
$\trianglelefteq_\lambda$-maximal.

What about reduced powers for $\lambda$-regular filter $D$ on
$\lambda$?   By \cite{Sh:17}, $M^\lambda/D$ is
$\lambda^+$-saturated for every $T$ (of cardinality $\le \lambda$) 
\Iff \, $D$ is $\lambda^+$-good
and $\cP(\lambda)/D$ is a $\lambda^+$-saturated Boolean algebra.  
Parallel results hold when we replace $\lambda^+$-saturated by
$(\lambda^+,\Sigma_{1+n}(\bbL_{\tau(T)}))$-saturated.  We 
shall concentrate on ($\lambda^+$, atomic)-saturated and 
introduce the
related partial
order $\trianglelefteq^{\rp}_\lambda$, see definitions below.

Concerning ultrapowers, lately, Malliaris-Shelah \cite{Sh:998} 
proved  that 
a regular ultra-filter $ D $ on  
a cardinal $ \lambda $ is $ \lambda ^+$-good
iff for any linear order $ M $ we have 
$ M^\lambda /D$ has no symmetric pre-cut with 
cofinality $ \le \lambda $. This was proved together with
the 
theorem $ \mathfrak{p} = \mathfrak{t} $
and 
``for a f.o. complete countable $T$,  being $\SOP_2$ suffices for $\triangleleft_\lambda$-maximality".
In a later work \cite{Sh:1051},  it is 
proved that 
at least for a relative
$\triangleleft^*_\lambda$ (see \cite{Sh:500}) this is ``iff" assuming
a case of G.C.H., relying also on works with Dzamonja \cite{Sh:692},
and with Usvyastov \cite{Sh:844}.  Part of
the proof is axiomatized by Malliaris-Shelah \cite{Sh:1070}.

Note also that \cite{Sh:1019} deals with saturation but only for
ultrapowers by $\theta$-complete ultra-filters for 
$\theta$ a compact cardinal; and also with $\omega$-ultra-limits.

Now what do we accomplish here?

First, in \S1 we axiomatize the proof of \cite[Ch.VI,2.6]{Sh:a},
i.e. we define when $\mathbf r = (M,\Delta)$ is
a so-called 
$\RSP$ and for it to prove
that the relevant model $N_{\mathbf r}$ is $(\min\{\gp_{\mathbf
  r},\gt_{\mathbf r}\},\Delta)$-saturated.  Second, in \S2 we prove, of
course, that \cite[Ch.VI,2.6]{Sh:a} follows, but also we show that the axiomatization of $\RSP$ is by Horn sentences.  Hence we can apply it to reduced
powers.  So $T$ is
$\trianglelefteq^{\rp}_\lambda$-maximal if $T=\Th(\bbQ,<)$ and moreover  for every $T$ having the $\SOP_3$; lastly, we comment on models of Peano Arithmetic.

In \S3 we try to sort out when for models of $T$ we get the relevant
atomic saturation.

Can we generalize also results \cite{Sh:998} to reduced powers?  

The main result of \S4 says that no.  We also sort out the parallel of goodness, excellence, and morality for filters 
and atomic saturation for reduced powers. 
In a hopeful continuation \cite{Sh:F1959}, we shall 
try to sort out the order $ \trianglelefteq ^{\rp}_\lambda $,  and in particular  consider non-maximality and 
parallel statements for infinitary logics (see \cite{Sh:1019}).

The reader can ignore Boolean ultra-powers (that is \ref{z9}, \ref{z10}, \ref{z12}   for sections 1,2 and can in first reading deal only with first-order logic (so $ \theta = {\aleph_0} $,  and the assumptions concerning the completeness of filer disappear. We thank the referee for many helpful comments.

Note that by \ref{c31}
\begin{conclusion}
\label{x0}
If $(T,\Delta)$ has the $\SOP_3$, \then\, it is
$\trianglelefteq^{\rp}_\lambda$-maximal. 
\end{conclusion}

\begin{question}
\label{x1}
Do we have:  if $D$ is $(\lambda_2,T)$-good and regular \then \, $D$ is
$(\lambda_1,T)$-good when $\lambda_1 < \lambda_2$ (or more)?
\end{question}

\subsection {Further Questions}\

\begin{convention}\label{x2}
    1) Let $T$ be a theory with elimination of quantifiers if not said otherwise.  Let $\Mod_T$ be the class of models of $T$.
    
    2) The main case is for $T$  a countable complete first-order theory with elimination of quantifiers, moreover, with every formula equivalent to an atomic one.
\end{convention}

So it is natural to ask

\begin{conjecture}
\label{x5}
The pair $(T,\Delta)$ is $\trianglelefteq_{\rp}$-maximal iff $(T,\Delta)$ 
has the $\SOP_3$.
\end{conjecture}

So which $T$ (with elimination of quantifiers) are maximal under
$\triangleleft^{\rp}_\lambda$?  That is, when for every regular filter
$D$ on $\lambda,M^\lambda/D$ is ($\lambda^+$, atomic)-saturated iff
$D$ is $\lambda^+$-good?  Is $T_{\feq}$ maximal?
(see \cite{Sh:457}, it is a proto-typical non-simple $ T $, 
but see more in \cite{Sh:1164})
As we have not
proved this even for ultra-filters, the reasonable hope is that it will
be easier to show non-maximality for $\triangleleft^{\rp}_\lambda$.  
Also in light of \cite{Sh:1030} for simple theories we like to
prove non-maximality with no large cardinals.  We may hope to use just
$\NSOP_2$, but still it would not settle the problem of characterizing
the maximal ones as e.g. $\SOP_2 \equiv \SOP_3$ is open for such
$T$; for  pairs $(T,\varphi(\bar x,\bar y))$ they are different.

Note that for first order $T$, it makes sense to use $\mu^+$-saturated
models and $D$ is $\mu^+$-complete.

Also, the ``$T$ stable" case should be resolved.

\begin{conjecture}
\label{b41}
$M^\lambda /D$ is ($\aleph^\lambda_0/D$, atomic)-saturated \when \,:

\begin{enumerate}
\item[$(a)$]  $T$ a theory as in \ref{x2}(2),
\sn
\item[$(b)$]  $T$ is stable without the $\fcp,$ 

\item[$(c)$]  $D$ is a 
regular filter on $\lambda$.
\end{enumerate}
\end{conjecture}

\begin{remark}
Maybe given a $1-\varphi$-type $p \subseteq \{\varphi(x,\bar a):
\bar a \in{}^m(M^I/D)\}$ of cardinality $\le \lambda$ in $M^I/D$, we
try just to find a dense set of $A \in D^+$ such that in $M^I/(D+A)$ the
$1-\varphi$-type is realized.  Then continue; opaque.
\end{remark}

\subsection {Preliminaries} \

\begin{notation}
\label{z0}
1) $ T $ dnote a f.o. theory, usually complete.

2) Let $ \tau $ denote a vocabulary,
 $ \tau _T= \tau (T)$ denote the vocabulary of the 
 theory $ T $

3) We use $ M,N$ to denote models, $ \tau_M= \tau(M)$
 is the vocabulary of $ M $ and $ P^M, F^M$ denote the 
 interpretation of $ P, F $ respectively.

 4) let $\mathbb{L} (\tau) $ denote the f.o. 
 language for the vocabulary $ \tau $.

 5) We allow function symbol $ F \in \tau $ to be interprated
 in a $ \tau $-model $ M $ as a partial
 function, but with domain $ P^M_F$, with 
 $ P_F \in \tau $ a predicate with the same
 arity.
\end{notation}

\begin{notation}\label{z1}
1) Let $\gB$ denote a Boolean algebra, $\comp(\gB)$ its completion, $\gB^+
= \gB \backslash \{0_\gB\},\uf(\gB)$ the set of ultra-filters on $\gB,
\fil(\gB)$ the set of filters on $\gB$.  For $\mathbf a \in \mathfrak{B} $ let
$\mathbf a^{\iif(\true)} = \mathbf a^{\iif(1)}$ be $\mathbf a$ and let
$\mathbf a^{\iif(\false)} = \mathbf a^{\iif(0)}$ be $1_{\gB} - \mathbf a$.

1A) Let $ \mathfrak{B}_1 \lessdot  \mathfrak{B} _2 $ mean that $ \mathfrak{B} _1 $ is a sub-algebra of $ \mathfrak{B}_2$, and moreover, a complete one, which means that every maximal anti-chain  of $ \mathfrak{B}_1 $ is a maximal anti-chain of $ \mathfrak{B}_2 $.
 
2) For a model $M$ let $\tau_M = \tau(M)$ be its vocabulary.

3) For a filter $ D$ on a set $ I $ let $D^+ = 
\{B \subseteq I: I \setminus B \notin D \}. $ 
\end{notation}

Now about cuts (they are
closed to but 
different than gaps, see
\cite{Sh:1069}).

\begin{definition}
\label{z2}
1) For a partial order $\cT = (\cT,\le_{\cT})$, we say $(C_1,C_2)$ is
pre-cut \when \, (but we may in this paper omit the {`}{`}pre"):

\begin{enumerate}
    \item[$(a)$]  $C_1 \cup C_2$ is a subset of $\cT$ linearly ordered by $\le_{\cT},$
    \sn
    \item[$(b)$]  if $a_1 \in C_1,a_2 \in C_2$ then $a_1 \le_{\cT} a_2,$
    \sn
    \item[$(c)$]  for no $c \in \cT$ do we have $a_1 \in C_1 \Rightarrow
    a_1 \le_{\cT} c$ and $a_2 \in C_2 \Rightarrow c \le_{\cT} a_2$.
\end{enumerate}

2) Above we say $(C_1,C_2)$ is a $(\kappa_1,\kappa_2)$-pre-cut \when \, in
addition:

\begin{enumerate}
    \item[$(d)$]  $C_1$ has cofinality $\kappa_1,$
    \sn
    \item[$(e)$]  $C^*_2$, the inverse of $C_2$, has cofinality $\kappa_2,$
    \sn
    \item[$(f)$]  so $\kappa_1,\kappa_2$ are regular
    cardinals,  
    (here 
    we ignore the case 0,1 if not said otherwise).  
\end{enumerate}

2A) Above we call $ \kappa_1, \kappa_2$ the cofinalities
of the pre-cut $ (C_1, C_2 ) $. We say that 
the pre-cut is symmetric when $ \kappa_1=\kappa_2$ 
and then we may say $ \kappa_1$ is its cofinality.

3) We may replace $C_\ell$ by a sequence $\bar a_\ell$, if not said
otherwise such that $\bar a_1$ is $\le_{\cT}$-increasing and $\bar a_2$ is 
$\le_{\cT}$-decreasing.

4) We say $(C_1,C_2)$ is a $(\kappa_1,\kappa_2)$-linear-cut of $\cT$ \when
\, it is a $(\kappa_1,\kappa_2)$-pre-cut and $C_1 \cup C_2$ is
downward closed, so natural for $\cT$ a tree.

5) We say $(C_1,C_2)$ is a weak pre-cut \when \, (b),(c) of part (1) holds.
%
\end{definition}

\begin{remark}
\label{z3}
1)
If $\cT$ is a (model-theoretic) tree, $\kappa_2 > 0$ and $(C_1,C_2)$
is a $(\kappa_1,\kappa_2)$-pre-cut \then \, it induces one and only
one $(\kappa_1,\kappa_2)$-linear-cut 
$(C'_1,C'_2)$, i.e. one satisfying $C_1 \subseteq
C'_1,C_2 \subseteq C'_2$ such that $C_1 \cup C_2$ is cofinal in $C'_1
\cup C'_2$.

2) In \ref{z6} below, if $ L
= 
\mathbb{L} (\tau ) $ 
then $ \theta = {\aleph_0} , \sigma =1$ suffice,
but not so in more general cases.
\end{remark}

\begin{definition}
\label{z6}
1) We say $M$ is fully $(\lambda,\theta,\sigma,L)$-saturated 
(may omit the fully); where $L \subseteq 
\cL(\tau_M)$ and $\cL$ is a logic; we may write 
$\cL$ if $L = \cL(\tau_M)$, \when \,: 

\begin{enumerate}
\item[$\bullet$]  if $\Gamma$ is a set of $< \lambda$ formulas from $L$ 
with parameters from $M$ with $< 1 + \sigma$ free variables, and
$\Gamma$ is $(< \theta)$-satisfiable in $M$, \then \,
$\Gamma$ is realized in $M$.
\end{enumerate}

2) We say ``locally" when using one $\varphi = \varphi(\bar x,\bar y) \in
\cL$ with $\ell g(\bar x) < 1 + \sigma$, i.e. all members of 
$\Gamma$ have the form\footnote{In \cite{Sh:1019} we use a $L
  \subseteq \bbL_{\theta,\theta},\theta$ a compact cardinal and if
  $\sigma > \theta$ we use a slightly different version of the
  definition of local and of the default value
 of $\sigma$ was
  $\theta$.}  $\varphi(\bar x,\bar b)$.

3) Saying ``locally/fully $(\lambda,\cL)$-saturated" the default
values (i.e. we may omit) of $\sigma$ is $\sigma = \theta$, 
of $(\sigma,\theta)$ is $\theta = \aleph_0 \wedge \sigma
= \aleph_0$ and of $\cL$ is $\bbL$ (first order logic) and of $L$ is $\cL$.
Omitting 
$\lambda$ means 
$\lambda = \|M\|$.

4) If $\varphi(\bar x,\bar{ y}) \in\cL(\tau_M)$ and $\bar a \in {}^{\ell
  g(\bar y)}M$ then $\varphi(M,\bar a) := \{\bar b \in {}^{\ell g(\bar
  x)}M:M \models \varphi [\bar b,\bar a]\}$.

5) Let $\bar x_{[u]} = \langle x_s:s \in u\rangle$.
\end{definition}

In \ref{z9}. \ref{z10}, \ref{z12} we shall deal with complete Boolean algebras and ultrapowers, and then we define 
an order between theories. 

\begin{definition}
\label{z9}
Assume we are given 
a Boolean algebra $\gB$ usually complete and a model or a set $M$
and $D$ a filter on $\comp(\gB)$, the completion of $\gB$.

1) Let $M^{\gB}$ be the set of partial functions $f$ from
 $\gB^+ 
 $ into $M$ such that 
for some maximal anti-chain $\langle
   a_i:i < i(*)\rangle$ of $\gB,\dom(f)$ includes $\{a_i:i < i(*)\}$
 and is included in\footnote{for the $D_\ell \in \uf(\gB_\ell)$ 
ultra-product, \wilog \, $\gB$ is complete, then
\wilog \, $f \rest \{a_i:i < i(*)\}$ is one to one.  But in general we
allow $a_i = 0_{\gB}$, those are redundant but natural in
\ref{z9}(3).} $\{a \in \gB^+
:(\exists i)(a \le a_i)\}$ and $f$ is a function into $M$ and $f \rest \{a \in
  \dom(f):a \le a_i\}$ is constant for each $i$.

1A) Naturally for $f_1,f_2 \in M^{\gB}$ we say $f_1,f_2$ are
$D$-equivalent, or $f_1 = f_2 \mod D$ \when \, for some $b \in D$ we have
$a_1 \in \dom(f_1) \wedge a_2 \in \dom(f_2) \wedge a_1 \cap a_2 \cap b
 > 0_{\gB} \Rightarrow f_1(a_1) = f_2(a_2)$.

2) We define $M^{\gB}/D$ naturally,
as well as $\TV_M (\varphi(f_0,\dotsc,f_{n-1})) \in \comp(\gB)$ \when \,
$\varphi(x_0,\dotsc,x_{n-1}) \in \bbL(\tau_M)$ and $f_0,\dotsc,f_{n-1} \in
 M^{\gB}$ where: 
 \begin{enumerate} 
 \item[(a)] $\TV$ stands for truth value, 
 \item[(b)] $ \TV_M  (\varphi(f_0,\dotsc,f_{n-1})) $ is equal to $$ 
 \sup \{a \in \mathfrak{B} ^+ : a \cap \bigcap _{{\ell} < n }
 \Dom( f_ {\ell} ) :
 M \models  (\varphi(f_0(a),\dotsc,f_{n-1} (a)))\}, $$ 
 \item[(c)] $ M $ is defined by letting, for $ \varphi $ 
 an atomic formulas 
 
 $M^{\gB}/D \models
\varphi[f_0/D,\dotsc,f_{n-1}/D]$ iff 
$\TV_M(\varphi(f_0,\dotsc,f_{n-1})) \in D$. 
\end{enumerate}

2A) Abusing notation, not only $M^{\gB_1} \subseteq M^{\gB_2}$ but
$M^{\gB_1}/D_1 \subseteq M^{\gB_2}/D_2$ when $\gB_1 \lessdot
\gB_2,D_\ell \in \fil(\gB_\ell)$ for $\ell=1,2$ and $D_1 = \gB_1 \cap
D_2$.  Also $[f_1,f_2 \in M^{\gB_1} \Rightarrow f_1 = f_2 \mod D_1
\leftrightarrow f_1 = f_2 \mod D_2]$.  So for 
$f \in M^{\gB_1}$ we identify $f/D_1$ and $f/D_2$.

3) For complete $\gB$, we say $\langle a_n:n < \omega\rangle$ represents $f \in
\bbN^{\gB}$ \when \, $\langle a_n:n < \omega\rangle$ is a maximal
anti-chain of $\gB$ (so $a_n = 0_{\gB}$ is allowed) 
and for some $f' \in \bbN^{\gB}$ which is 
$D$-equivalent to $f$ (see \ref{z9}(1A)) we have $f'(a_n)=n$.

4) We say $\langle (a_n,k_n):n < \omega\rangle$ represents  
$f \in
\bbN^{\gB}$ \when \,:

\begin{enumerate}
\item[$(a)$]  the $k_n$ are natural numbers with no repetition,
\sn
\item[$(b)$]  $\langle a_n:n < \omega \rangle$ is a maximal anti-chain,
\sn
\item[$(c)$]  $f(a_n) = k_n$.
\end{enumerate}

5) If $\cI$ is a maximal anti-chain of $\gB$ and $\bar M = \langle
M_a:a \in \cI\rangle$ is a sequence of $\tau$-models, then we define
$\bar M^{\gB}$ as 
the set of partial functions $f$ from $\gB^+$ to
$\cup\{M_a:a \in \cI\}$ such that for some maximal anti-chain $\langle
a_i:i < i(*)\rangle$ of $\gB$ refining $\cI$ (i.e. $(\forall i <
i(*))(\exists b \in \cI)(a_i \le_{\gB} b))$ we have:

\begin{enumerate}
\item[$(a)$]  $\{a_i:i <i(*)\} \subseteq \dom(f) \subseteq \{b \in
  \gB^+:b \le_{\gB} a_i$ for some $i <i(*)  \} $
\sn
\item[$(b)$]  if $a \in \dom(f)$ and $a \le a_i$ then $f(a) = f(a_i)$
\sn
\item[$(c)$]  if $a_i \le_{\gB} b,b \in \cI$ then $f(a_i) \in M_b$.
\end{enumerate}

6) For $\bar M,\gB,\cI$ as above and a filter $D$ on $\gB$ we define
$\bar M^{\gB}/D$ as in part (2) replacing $M^{\gB}$ there by $\bar
M^{\gB}$ here, see part (7).

7) For $\bar M,\gB,\cI$ as above, $\varphi = \varphi(\bar x) =
\varphi(x_0,\dotsc,x_{n-1}) \in \bbL(\tau_M)$ and $\bar f = \langle
f_\ell:\ell < n\rangle$ where $f_0,\dotsc,f_{n-1} \in \bar M^{\gB}$,
let $\TV(\varphi[\bar f]) = \TV(\varphi[\bar f],\bar M^{\gB})$ be
$\sup\{a \in \gB^+$: if $\ell < n$ then $a \in \dom(f_\ell)$ and $a
\le b \in \cI$ then $M_b \models \varphi[f_0(b),\dotsc,f_{n-1}(b)]
\} $.

8) We say $\gB$ is 
$  (< \sigma  ) $-distributive when it is $ \theta $-distributive
for every $ \theta < \sigma $, where

8A) $ \mathfrak{B} $ 
is 
$\theta$-distributive \when \,: if
for  
$\alpha <
\theta,\cI_\alpha $ is a maximal antichain of $\gB$ 
\then \, there is a maximal antichain of $\gB$ refining every
$\cI_\alpha(\alpha < \theta )$; 
this holds, e.g. when 
$\gB = \cP(\lambda)$ or just
there is a dense $Y \subseteq \gB^+$ closed under intersection of $
\theta$. 
\end{definition}

\begin{definition}
\label{z10}
1) 
Let $\gB$ be a complete Boolean algebra and $D$ a filter on $\gB$.  We
say that $D$ is $(\mu,\theta)$-regular \when \, for some 
$(\bar{\mathbf c},\bar{ u } )$ we have:

\begin{enumerate}
\item[$(a)$]  $\bar{\mathbf c} = \langle \mathbf c_\alpha:\alpha <
  \alpha_*\rangle$ is a maximal anti-chain of $\gB,$
\sn
\item[$(b)$]  $\bar u = \langle u_\alpha:\alpha < \alpha_*\rangle$
  with $u_\alpha \in [\mu]^{< \theta},$
\sn
\item[$(c)$]  if $i < \mu$ then $\sup\{\mathbf c_\alpha:\alpha$ satisfies $i
  \in u_\alpha\} \in D$.
\end{enumerate}

2) A filter $ D $ is called $ \lambda $-regular when 
it is $ (\lambda , {\aleph_0}  ) $-regular; the filter 
$ D$
on a set $ I $ (that is the Boolean algebra $ {\mathscr P} (I)$)
is called regular when it is a filter on a set $ I $ 
and it is $ |I|$-regular.
\end{definition}

\begin{claim}
\label{z12}
Assume $\gB$ is a complete Boolean
algebra 
which is 
$(< \lambda ) $-distributive and 
$D$ a filter on $\gB$ and $\theta = \cf(\theta) \le \lambda $.

1) Assume 
$ D$  is a $\theta $-complete 
ultra-filter. 
The parallel of {\L}o\'s 
theorem holds for
$\bbL_{\theta,\theta}$ and if $D$ is $\lambda$-complete even for
$\bbL_{\lambda,\theta}$ which means: if $\bar M = \langle M_b:b \in
\cI\rangle$ is a sequence of $\tau$-models, $\cI$ is a maximal
antichain of the complete Boolean algebra $\gB$ and $\varepsilon <
\theta,\varphi = \varphi(\bar x_{[\varepsilon]}) \in
\bbL_{\lambda,\theta}(\tau)$ and $f_\zeta \in \bar M^{\gB}$ for $\zeta
< \varepsilon$ then $M^\mathfrak{B} /D \models ``\varphi[\langle f_\zeta/D:\zeta <
\varepsilon\rangle]"$ iff $\TV_M (\varphi[\langle f_\zeta/D:\zeta <
\varepsilon\rangle])$ belongs to $D$.

2) If  in addition 
$D$ is $(\lambda, \theta ) $-regular and $M,N$ are
$\bbL_{\theta,\theta}$-equivalent 
\then \,  $M^{\gB}/D,N^{\gB}/D$ are
$\bbL_{\lambda^+,\theta}$-equivalent. 
\end{claim}


\begin{definition}
\label{z16}
1) Assume $\Delta_\ell$ is a of set atomic formulas in $\bbL(\tau(T_\ell))$.
Then we say $(T_1,\Delta_1) \trianglelefteq^{\rp}_{\lambda,\theta} 
(T_2,\Delta_2)$ \when \,: if $D$ is
  a $(\lambda, \theta ) $-regular filter on $\lambda$ and $M_\ell$ is a
  $\lambda^+$-saturated model of $T_\ell$ for $\ell=1,2$ and
  $M^\lambda_2/D$ is $(\lambda^+,\theta,\Delta_2)$-saturated \then \,
  $M^\lambda_1/D$ is $(\lambda^+,\theta,\Delta_1)$-saturated.

2) For general $\Delta_1,\Delta_2$ we define $(T_1,\Delta_1)
\trianglelefteq^{\rp}_{\lambda,\theta} (T_2,\Delta_2)$ as meaning
$(T^+_1,\Delta^+_1) \trianglelefteq^{\rp}_{\lambda,\theta}
(T^+_2,\Delta^+_2)$ where (as Morley \cite{Mo65} does):

\begin{enumerate}
\item[$\bullet$]  $T^+_\ell = T_\ell \cup \{(\forall \bar
  x)(\varphi(\bar x) \equiv P_\varphi(\bar x)):\varphi(\bar x) \in
  \Delta_\ell\}$ with $\langle P^\ell_\varphi:\varphi \in
  \Delta_\ell\rangle$ new pairwise distinct predicates with suitable
  number of places
\sn
\item[$\bullet$]  $\Delta^+_\ell = \{P^\ell_\varphi(\bar x_\varphi):
\varphi \in \Delta_\ell\}$.
\end{enumerate}

3) In (2), $T_1 \trianglelefteq^{\rp}_{\lambda,\theta} T_2$ means
$\Delta_\ell =$ the set of atomic
$\bbL_{\theta,\theta}(\tau_{T_\ell})$-formulas. 
\end{definition}

\begin{observation}
\label{z19}
Assume $\Delta \subseteq \bbL(\tau_T)$ is closed under $\exists$ and
$\wedge$.  A model $M$ of $T$ is $(\mu^+,\mu^+,\Delta)$-saturated \Iff
\, it is $(\mu^+,1,\Delta)$-saturated.
\end{observation}

\begin{question}   
1)
 Under $ \trianglelefteq _{\rp}$ characterize 
the minimal/maximal pairs $(T, \Delta ). $

2) What about the   
parallel of $ \trianglelefteq^{**}$   
(see \cite{Sh:457}, \cite{Sh:1051})?
\end{question}  
\newpage

\section {Axiomatizing \cite[Ch.VI,2.6]{Sh:c}} \label{1}

Note that while the notation $\gt(\cT)$ is obviously natural the
notation $\gp(\cT)$ is really justified 
just
by the results here.

\begin{definition}  
\label{h2}
1) For  a partial order $\cT = (\cT,\le_{\cT})$ 
let $\gp_{\cT} = \gp(\cT)$ be $\min\{\kappa_1 +
\kappa_2:(\kappa_1,\kappa_2) \in \cC_{\cT}\}$ and $\gp_\theta(\cT) =
\min\{\kappa_1 + \kappa_2:(\kappa_1,\kappa_2) \in \cC_{\cT,\theta}\}$;
where:

2)
$ \cC_\theta (\cT) = 
\{(\kappa_1,\kappa_2)$: the partial order $\cT$
has a $(\kappa_1,\kappa_2)$-cut
 and $\kappa_1 \ge \theta,\kappa_2 \ge \aleph_0\}$. 
 If $ \theta = {\aleph_0} $ then we may omit $ \theta $, 
 (yes, when $ \theta > {\aleph_0} $ this is not symmetric).

3) For a partial order $\cT$ let $\gt_{\cT} = \gt(\cT)$ be the minimal $\kappa \ge \aleph_0$ such that there is a $<_{\cT}$-increasing sequence 
of length $\kappa$ with no $<_{\cT}$-upper bound.

4) Let $\gp^*_{\cT} = \gp^*(\cT)$ be $\min\{\gt_{\cT},\gp_{\cT}\}$.

5) $\gp_{\theta-\sym} (\cT) =
\min\{\kappa:(\kappa,\kappa) \in \cC_\theta (\cT)\}$ and if $ \theta = {\aleph_0} $ we may write $\gp^*_{\sym}(\cT).$ 

6) In Definition \ref{h4} below let $\gt_{\mathbf r} = \gt_{\cT_{\mathbf r}},
\gp_{\mathbf r} = \gp_{\theta_{\mathbf r}}
(\cT_{\mathbf r})$.
\end{definition}

\begin{definition}  
\label{h4}
1) For $\iota = 1,2$ 
(the difference is only in closed (i)), 
we say $\mathbf r$ or $(M,\Delta)$ is a
$(\theta,\iota)$-realization\footnote{When $P$ and $\tau_N$ (hence
  $N$) are understood from the context we may omit them.} spectrum problem, in short
$(\theta,\iota)-\RSP$ or $(\theta,\iota)-1$-$\RSP$ \when \, $\mathbf r$ 
consists of (if $\iota =2$ we may omit it, similarly if $\theta =
\aleph_0$; we may omit $\Delta$ and write $M$ when $\Delta$ is the set of
atomic formulas in $\bbL(\tau_{N_M})$, see below, so $M$ below $=
M_{\mathbf r}$, etc.):

\begin{enumerate}
\item[$(a)$]  $M$ a model,
\sn
\item[$(b)$]  for the relations $\cT = \cT^M,\le_{\cT} = \le^M_{\cT}$ of $M$
  (i.e. $\cT,\le_{\cT}$ are predicates from $\tau_M$) we have
$\cT = (\cT,\le_{\cT})$ a partial order (so definable in $M$) with
  root $c^M = \rt(\cT)$, so $c \in \tau_M$ is an individual constant
  and $t \in \cT \Rightarrow \rt(\cT) \le_{\cT} t$; as in other cases
we may write $\cT_{\mathbf r},\le_{\mathbf r}$ for $\cT,\le_{\cT}$; we do
not require $\cT$ to be a tree; but do require $t \in \cT
\Rightarrow t \le_{\cT} t,$

\sn
\item[$(c)$]  a model $N = N_{\mathbf r} = N_M$ 
with universe $P^M,\tau(N) \subseteq \tau(M)$ such that
\sn
\begin{enumerate}
\item[$\bullet$]  $Q \in \tau_N \Rightarrow Q^M = Q^N,$
\sn
\item[$\bullet$]  $F \in \tau_N \Rightarrow F^N = F^M$, (we
  understand $F^M, F^N$ to be partial functions), 
\newline
so every $\varphi \in \bbL(\tau_N)$ can be
interpreted as $\varphi^{[*]} \in \bbL(\tau_M)$, all variables varying
on $P$ (include quantification); we may forget the 
    superscript  
$[*]$.
\end{enumerate}
\sn
\item[$(d)$]  the cardinal $\theta$ and  $\Delta \subseteq \{\varphi:\varphi = \varphi(x,\bar y) \in \bbL_{\theta,\theta}(\tau_N)\}$ which is closed under
conjunctions meaning: if $\varphi_\ell(x,\bar y_\ell) \in 
\Delta$ for $\ell=1,2$ then
 $\varphi(x,\bar y'_1,\bar y''_2) = \varphi_1(x,\bar y'_1)
\wedge \varphi_2(\bar x,\bar y''_2) \in \Delta,$
\sn
\item[$(e)$]  $R^M \subseteq |N| \times \cT^M$ so a two-place relation;
and let $R^M_t = \{b:b R^M t\}$ for $t \in \cT^M,$
\sn
\item[$(f)$]  $|N| \times \{\rt_{\cT}\} \subseteq R^M$,
  i.e. $R^M_{\rt(\cT)} = |N|,$
\sn
\item[$(g)$]  if $s \le_{\cT} t$ then $a \in N \wedge aRt \Rightarrow
  aRs$, i.e. $R^M_s \supseteq R^M_t,$
\sn
\item[$(h)$]  $t \in \cT \Rightarrow R^M_t \ne \emptyset,$
\sn
\item[$(i)$]  if $s \in \cT,\varphi(x,\bar a) \in \Delta(N) :=
  \{\varphi(x,\bar a)):\varphi(x,\bar y) \in \Delta$ and $\bar a \in {}^{\ell
g(\bar y)}N\}$ and for some $b \in R^M_s,N \models \varphi[b,\bar a]$ 
\then \, there is $t \in \cT$ such that $s \le_{\cT} t$ and
$R^M_t = \{b \in R^M_s:N \models \varphi[b,\bar a]\},$
\sn
\item[$(i)^+$]  if $\iota=1$ like clause $(i)$ 
but\footnote{We may not add a function, maybe it matters  when we try to build $\mathbf r$ with $\Th(M_{\mathbf
  r})$  nice first order.} moreover $t = F^M_{\varphi,1}(s,\bar a)$
where $F^M_{\varphi,1}:\cT_{\mathbf r} \times {}^{\ell g(\bar y)}(P^M)
  \rightarrow \cT_{\mathbf r},$
\sn
\item[$(j)$]  if $t \in \cT_{\mathbf r}$ and $\varphi(x,\bar a) \in
  \Delta(N)$ and $\varphi(N,\bar a) \ne \emptyset$ \then
\sn
\begin{enumerate}
\item[$(\alpha)$]  $s = F^M_{\varphi,2}(t,\bar a)$ is such that $R^M_s
  \cap \varphi(N,\bar a) \ne \emptyset$ and $s \le_{\cT} t,$
\sn
\item[$(\beta)$]  if $s = F^M_{\varphi,2}(t,\bar a),s_1 \le_{\cT} t$
  and $R^M_{s_1} \cap \varphi(N,\bar a) \ne \emptyset$ then $s_1
  \le_{\cT} s.$
\end{enumerate}
\sn
\item[$(k)$]  if $\theta > \aleph_0$ then in
$(\cT,\le_{\cT})$ any increasing chain of length $< \theta$
which has an upper bound has a $\le_{\cT}$-lub.
\end{enumerate}
\end{definition}

\begin{remark}
\label{h5}
We may consider adding: $S^M$ a being successor, (but this is not
Horn), i.e.:

\begin{enumerate}
\item[$(l)$] if $\iota =1$ we also have $S^M $ 
is the set of pairs 
$ (a,b)$ such that $ b $   
is a
  $\le_{\cT}$-successor of $a$  
    which means:  
\sn
\begin{enumerate}
\item[$(\alpha)$]  if $a \le b \wedge a \ne b$ then for some
  $c,S(a,c) \wedge c \le b,$
\sn
\item[$(\beta)$]  if $b \in \cT \backslash \{\rt_{\cT}\}$ then for
  some unique $a$ we have $S^M(a,b),$
\sn
\item[$(\gamma)$]  $S(a,b) \Rightarrow a \le b,$
\sn
\item[$(\delta)$]  $S(a,b_1) \wedge S(a,b_2) \wedge b_1 \ne b_2
  \Rightarrow \neg(b_1 \le b_2),$
\sn
\item[$(\varepsilon)$]  in clause (j) we can add $S^M(s,t)$.
\end{enumerate}
\end{enumerate}
\end{remark}

\begin{remark}  
\label{h5d}
Presently, it may be that $a \le_{\cT} b \le_{\cT} a$ but $a \ne b$.  Not a
disaster to forbid but no reason.
\end{remark}

How does this axiomatize realizations of types?
\begin{cd}  
\label{h6}
Let $\iota = \{1,2\},\theta$ is $\aleph_0$ or just a regular cardinal.

1) For any model $N$ and $\Delta \subseteq \{\varphi:
\varphi =  
   \varphi(x,\bar y) \in \bbL_{\theta,\theta}(\tau_T)\}$ closed under
   conjunctions of $< \theta$, the canonical $(\theta,\iota)-\RSP,
   \mathbf r = \mathbf r^\theta_{N,\Delta}$ defined below is indeed a 
$\theta-\RSP$.

2) $\mathbf r = \mathbf r^\theta_{N,\Delta}$ (if $\theta = \aleph_0$ we may
omit it) is defined by:

\begin{enumerate}
\item[$(a)$] $\Delta_{\mathbf r} = \Delta,N_{\mathbf r} = N$ and
  $\theta_{\mathbf r} = \theta,$
\sn
\item[$(b)$]  $\cT_{\mathbf r} = \{\langle \varphi_\varepsilon(x,
\bar a_\varepsilon):\varepsilon < \zeta\rangle:\zeta < \theta$ 
and for every $\varepsilon < \zeta$ we have 
$\varphi_\varepsilon(x,\bar a_\varepsilon) \in \Delta(N)$ 
and $N \models (\exists x)(\bigwedge\limits_{\varepsilon <
  \zeta} \varphi_\varepsilon(x,\bar a_\varepsilon))\},$
\sn
\item[$(c)$]  $\le_{\mathbf r} =$ being the initial segment relation
  on $\cT_{\mathbf r},$
\sn
\item[$(d)$]  $M = M_{\mathbf r}$ is the model with universe 
$\cT_{\mathbf r} \cup |N|$; \wilog \, $\cT_{\mathbf r} \cap |N| =
\emptyset$, with the relations and functions of $N,\cT_r,\le_{\mathbf r}$ and,
\sn
\begin{enumerate}
\item[$\bullet$]  $P^M = |N|,$
\sn
\item[$\bullet$]  $c^M = \langle \rangle \in \cT_{\mathbf r},$
\sn
\item[$\bullet$]  $R^M = \{(b,t):  b 
\in N,t = \langle
  \varphi_{t,\varepsilon}(x,\bar a_{t,\varepsilon}):\varepsilon 
< \zeta_t\rangle \in \cT$ and $N \models \varphi_{t,\ell}
(b,\bar a_{t,\varepsilon})$ for every $\varepsilon < \zeta_t\},$
\sn
\item[$\bullet$]  $F^M_{\varphi,2}$ as in Definition \ref{h4}(j),
\sn
\item[$\bullet$]  if $\iota=1$ then $F^M_{\varphi,1}$ is as in
Definition \ref{h4}(i)$^+$.
\end{enumerate}
\end{enumerate}
\end{cd}

\begin{remark}
\label{h7}
If we adopt \ref{h5} it is natural to add:

\begin{enumerate}
\item[$(e)$]  for $\iota = 1,S^M =
  \{(\bar\varphi_1,\bar\varphi_2):\bar\varphi_2 = 
  \bar{ \varphi}_1    
  \char 94
  \langle \varphi(x,\bar a)\rangle \in \cT_{\mathbf r}$ 
for some $\varphi(x,\bar a) \in \Delta(N)\}$.
\end{enumerate}
\end{remark}

\begin{PROOF}{\ref{h6}}  
Obvious.
\end{PROOF}

\begin{mc}  
\label{h8}
1) Assume $\mathbf r$ is an $\RSP$.  If $\kappa = \min\{\gt_{\mathbf r},
\gp_{\mathbf r}\}$ \then \, the
model $N$ is $(\kappa,1,\Delta_{\mathbf r})$-saturated, i.e.

\begin{enumerate}
\item[$\oplus$]  if $p(x) \subseteq \Delta_{\mathbf r}(N_{\mathbf r})$ is
  finitely satisfiable in $N_{\mathbf r}$ (= is a type in $N_{\mathbf r}$) of
  cardinality $< \kappa$ \then\, $p$ is realized in $N_{\mathbf r}$.
\end{enumerate}

2) If $\theta > \aleph_0$ and $\mathbf r$ is a $\theta$-$\RSP$, 
\then \, $N_{\mathbf r}$ is
$(\kappa,1,\Delta_{\mathbf r})$-saturated where $\kappa =
\min\{\gt_{\mathbf r},\gp_{\mathbf r}\}$ recalling \ref{h2}(6),
i.e. $\gp_{\mathbf r} = \gp_{\cT_{\mathbf r},\theta}$.

3) If $\theta > \aleph_0,\mathbf r$ is a $\theta$-$\RSP$ satisfying $(k)^+$
below \then \, $N_{\mathbf r}$ is $(\gt_{\mathbf r},1,
\Delta_{\mathbf r})$-saturated \when \,:

\begin{enumerate}
\item[$(k)^+$]  in $(\cT,\le_{\cT})$ any increasing chain which has an
upper bound, has a $\le_{\cT}$-lub.
\end{enumerate}
\end{mc}

\begin{PROOF}{\ref{h8}}  
This is an abstract version of \cite[Ch.VI,2.6]{Sh:a} = \cite[Ch.VI,2.6]{Sh:c};
recall that \cite[Ch.VI,2.7]{Sh:a} translates trees to linear orders.

1) Let $N = N_{\mathbf r},\Delta = \Delta_{\mathbf r}$, etc.

Let $p$ be a $(\Delta,1)$-type in $N$ of cardinality $< \kappa$.
\Wilog \, $p$ is infinite and closed under conjunctions.

So let

\begin{enumerate}
\item[$(*)_1$]  $\alpha_* < \kappa,p = \{\varphi_\alpha(x,\bar
  a_\alpha):\alpha < \alpha_*\} \subseteq \Delta(N),p$ is finitely
  satisfiable in $N$.
\end{enumerate}

We shall try to choose $t_\alpha$ by induction on $\alpha \le \alpha_*$
such that:

\begin{enumerate}
\item[$(*)_2$]  $(a) \quad t_\alpha \in \cT$ and $\beta < \alpha
  \Rightarrow t_\beta \le_{\cT} t_\alpha,$
\sn
\item[${{}}$]  $(b) \quad$ if $\beta < \alpha_*$ then there
  is $b \in R^M_{t_\alpha}$ such that $N \models \varphi_\beta[b,\bar
  a_\beta],$
\sn
\item[${{}}$]  $(c) \quad$ if $\beta < \alpha$ then $b \in
  R^M_{t_\alpha} \Rightarrow N \models \varphi_\beta[b,\bar a_\alpha]$.
\end{enumerate}

If we succeed, this is enough because if $t = 
t_{\alpha_*}$ is well defined then
$R^M_t \ne \emptyset$ by Definition \ref{h4}(h) and any 
$b \in R^M_t$ realizes the type by $(*)_2(c)$ and Definition
\ref{h4}(h).   Why can we carry the definition?

\underline{Case 1}:  $\alpha = 0.$

Let $t_\alpha = \rt_{\cT}$, hence $R^M_{t_\alpha} = |N|$ by Definition
\ref{h4}(f).  Now clause (a) of $(*)_2$ holds as $t_\alpha \in 
\cT_{\mathbf r}$ and there is no $\beta <
\alpha$.  Also clause (b) of $(*)_2$ holds because $p$ is a type
and $R^M_{\rt(\cT)} = |N_{\mathbf r}|$ by Definition \ref{h4}(h).

Lastly, clause (c) of $(*)_2$ holds trivially.

\underline{Case 2}: $\alpha = \beta +1.$

If $\iota =1$ let $t = F^M_{\varphi_\beta,1}(t_\beta,\bar a_\beta)$ 
and see clause (i)$^+$ of Definition \ref{h4}.  If $\iota=2$ use
clause (i) of the definition recalling $p$ is closed under conjunctions.

\underline{Case 3}:  $\alpha$ a limit ordinal

As $\gt_{\cT_{\mathbf r}} \ge \kappa > \alpha_*$ by the claim's 
assumption (on $\gt_{\cT_{\mathbf r}}$, see Definition \ref{h2}(2)) 
necessarily there is $s \in \cT$ such that
$\beta < \alpha \Rightarrow t_\alpha \le_{\cT} s$.  We now try to choose
$s_i$ by induction on $i \le \alpha_*$ such that

\begin{enumerate}
\item[$(*)_{2.1}$]  $(a) \quad s_i \in \cT,$
\sn
\item[${{}}$]  $(b) \quad \beta < \alpha \Rightarrow t_\beta \le_{\cT} s_i,$
\sn
\item[${{}}$]  $(c) \quad j < i \Rightarrow s_i \le_{\cT} s_j,$
\sn
\item[${{}}$]  $(d) \quad$ if $i = j+1$ then $R^M_{s_i}$ is not
  disjoint to $\varphi_j(N,\bar a_j)$.
\end{enumerate}

If we succeed, then $s_{\alpha_*}$ satisfies all the demands on
$t_\alpha$ (e.g. $(*)_2(b)$ holds by Definition \ref{h4}(g) and
$(*)_{2.1}(d))$,  so we have just to carry the induction for $\alpha$.  
Now if $i=0$ clearly $s_0 = s$ it as required.  If $i
=j+1$ let $s_i = F^M_{\varphi,2}(s_j,\bar a_j)$, by Definition
\ref{h4}(j) it is as required.  For $i$ a limit ordinal use 
$\kappa \le \gp_{\cT}$ hence
carry the induction on $i$ so finish
case 3.  

So we succeed to carry the induction on $\alpha$ hence (as said after $(*)_2$) 
get the desired conclusion.

2) Similar, except concerning case 3.  Note that \wilog \, $\theta >
   \aleph_0$ by part (1).

\underline{Case 3A}:  $\alpha$ is a limit ordinal of cofinality $\ge \theta.$

As in the proof of part (1).

\underline{Case 3B}:  $\alpha$ is a limit ordinal of cofinality $< \theta.$

Again there is an upper bound $s$ of $\{t_\beta:\beta < \alpha\}$.  Now
by clause (k) of Definition \ref{h4}, \wilog \, $s$ is a $<_{\cT}$-lub
of $\{t_\beta:\beta < \alpha\}$.  So easily for every $i <
\alpha_*,F^N_{\varphi_i,2}(s,\bar a_i)$ is $\ge t_\beta$ for $\beta <
\alpha$ hence is equal to $s$, so $s_\alpha := s$ is as required.

3) Similarly.
\end{PROOF}

\begin{discussion}\label{h10}
1) What about ``$(\lambda^+,n,\Delta)$-saturation"?  We can repeat the
same analysis or we can change the models to code $n$-tuples.  More
generally, replacing $\varphi(\bar x_{[\varepsilon]},\bar y)$ 
by $\varphi(\langle F_\zeta(x):\zeta < \varepsilon\rangle,\bar y)$,
using $F_\zeta \in \tau_M$ (though not necessarily $F_\zeta \in
\tau_{N_{\mathbf r}}$), so we can allow infinite $\varepsilon$.

2) Hence the same is true for
$(\lambda^+,\aleph_0,\Delta)$-saturation, e.g.
$\lambda^+$-saturated by an assumption.
\end{discussion}

\newpage

\section{Applying the axiomatized frame} \label{2}

%
Consider a filter $ D $ on a set $ I $ and cardinals
$ \lambda \ge \mu $. We may ask for a model $ M $ 
of cardinality $ \ge \mu $, whether 
$ M^I/D $  is 
$ (\lambda ^+ , \rm{atomic})$-saturated, varying
$ M $.

We here apply \S1 to show that: when $ D $ is an ultra-filter, 
the model $ ({}^{ \omega > } \mu , \triangleleft )$ is the
hardest, this is \ref{c4}, We then (in \ref{c7}) show that 
\S1 has axiomatization 
which is a 
Horn theory. Hence we can prove results
like \ref{c4} below  for filters $ D $ (not just 
for ultra-filters),

\begin{conclusion}
\label{c4}
1) If $D$ is an ultra-filter on a set $I,N$ a model, $\mu = \|N\| +
\|\tau_N\|$ and $({}^{\omega >}\mu,\triangleleft)^I/D$ is
($\lambda^+$,atomic)-saturated \then \, $N^I/D$ is $\lambda^+$-saturated.

2) Instead of  ``$({}^{\omega >}\mu,\trianglelefteq)^I/D$ is ($\lambda^+,1$,
atomic)-saturated" we can demand  ``$J^I/D$ is
($\lambda^+,1$,atomic)-saturated" where $J$ is the linear order with
set of elements $\{-1,1\} \times {}^{\omega >}\mu$ ordered by
$(\iota_1,\eta_1) < (\iota_2,\eta_2)$ \Iff \, $\iota_1 < \iota_2$ or
$\iota_1 = -1 = \iota_2 \wedge \eta_1 <_{\lex} \eta_2$ or $\iota_1 =
-1 = \iota_2 \wedge \eta_2 <_{\lex} \eta_1$.
\end{conclusion}

\begin{PROOF}{\ref{c4}}
1) Let $N_1 = N$.  
As $D$ is an ultra-filter \wilog \, $\Th(N_1)$ has elimination of
quantifiers and even every formula is
equivalent to an atomic formula.  
Let $\Delta = \bbL(\tau_N)$, by \ref{h6}  $\mathbf
r_1 := \mathbf r_{N_1,\Delta}$ is an $\RSP$.  Let $N_2 = N^I_1/D $ 
and let $ M_1 =  
M_{\mathbf r_1},M_2  = M^I_1/D$ and let $\mathbf r_2$ be the 
$\RSP(M_2 $ , $ \Delta)$.  Clearly, $\mathbf r_2$ is an $\RSP$ as the demands in \ref{h4} are first order (see more in \ref{c7}).  

Now,

\begin{enumerate}
\item[$(*)_1$]  $\cT_{\mathbf r_1} \cong ({}^{\omega
    >}\mu,\triangleleft)$.
\end{enumerate}

[Why?  See \ref{h6}(2).]

\begin{enumerate}
\item[$(*)_2$]  $\cT_{\mathbf r_2} = (\cT_{\mathbf{r} _1})^I/D$
is ($\lambda^+$,atomic)-saturated.
\end{enumerate}

[Why?  By an assumption.]

\begin{enumerate}
\item[$(*)_3$]   $\gt(\cT_{\mathbf r_1}),\gp(\cT_{r_2}) \ge \lambda^+$.
\end{enumerate}

[Why?  Follows by $(*)_2$.]

Hence by \ref{h8}, $N_2$ is $(\lambda^+,1,1,\Delta)$-saturated which
means $N_2 = (N_1)^I/D$ is $\lambda^+$-saturated.

2) Easy (or see \cite[Ch.VI,2.7]{Sh:a} or see \cite{Sh:14}).
\end{PROOF}

To apply the criterion of the Main Claim \ref{h8} to reduced products
we need:
\begin{claim}
\label{c7}
If $\Delta$ is the set of conjunctions of atomic formulas (no
negation!) in $\bbL(\tau_0)$ and $\tau = \{\cT,\le_{\cT},R,P,c\} \cup
\{F_{\varphi,\ell}:\varphi \in \Delta$ and $\ell=2$ or 
$\ell=1$ if relevant$\} \cup \tau_0$,  
(disjoint union, recall $c$ is $\rt_{\cT}$), \then \, there is    a 
set $T$ of Horn
sentences from $\bbL(\tau)$ such that for every $\tau$-model $M$

\begin{enumerate}
\item[$\bullet$]  $(M,\Delta)$ is a $\RSP$ (i.e. 2-$\RSP$) \Iff \, $M
  \models T$.
\end{enumerate}
\end{claim}

\begin{PROOF}{\ref{c7}}
Consider Definition \ref{h4}.  For each clause we consider the
sentences expressing the demands  
there.

\underline{Clause (a)}:  Obvious.

\underline{Clause (b)}:  Clearly the following are Horn:

\begin{enumerate}
\item[$\bullet$]  $x \le_{\cT} y \rightarrow \cT(x)$ 
    and 
$ x \le_{\cT} 
y \rightarrow \cT(y),$
\sn
\item[$\bullet$]  $x \le_{\cT} y \wedge y \le_{\cT} z \rightarrow x
  \le_{\cT} z$,
\sn
\item[$\bullet$]  $\cT(\rt_{\cT})$ and $\cT(s)
\rightarrow \rt_{\cT} \le s,$
\sn
\item[$\bullet$]  $\cT(x) \rightarrow x \le_{\cT} x$.
\end{enumerate}

Note that $(\cT,\le_{\cT})$ being a tree is not a Horn sentence but
this  
is
not required.

\underline{Clause (c)}:

\begin{enumerate}
\item[$\bullet$]  $Q(x_0,\dotsc,x_{n(Q)-1}) \rightarrow
P(x_\ell)$ when $Q$ is an $n(Q)$-place predicate from $\tau(N)$ and
$\ell < n(Q)$; clearly it is Horn,
\sn
\item[$\bullet$]  for any $n$-place function symbol $F \in \tau_0$ the sentence:
$P(x_0) \wedge \ldots \wedge P(x_{n-1}) \rightarrow
P(F(x_0,\dotsc,x_{n-1}))$ and $y = F(x_0,\dotsc,x_{n-1}) \rightarrow
P(x_\ell)$.  
\end{enumerate}

\underline{Clause (d)}:  nothing to prove - 
see the present claim assumption on $\Delta$.

Recall  that for $ F \in \tau _N$, $ F $ stand for a partial 
function symbol with domain $ P_F$. 

\underline{Clause (e)}:  $yRs \rightarrow \cT(s),yRs 
\rightarrow P(y)$ are Horn.

\underline{Clause (f)}: $P(x) \rightarrow xR(\rt_{\cT})$ is Horn.

\underline{Clause (g)}: $s \le_{\cT} t \wedge xRt \rightarrow xRs$ is Horn.

\underline{Clause (h)}:  $(\forall t)(\exists x)(\cT(t) \rightarrow
xRt)$ is Horn.

\underline{Clause (i)}:  Let $\varphi(x,\bar y) \in \Delta$.  

First assume $\iota = 1$.  Note the following are Horn: for any
$\varphi(x,\bar y) \in \Delta,$

\begin{enumerate}
\item[$\bullet$]  $\cT(s) \wedge xRs \wedge \varphi(x,\bar y) \wedge
 \bigwedge\limits_{\ell < \ell g(\bar y)} P(y_\ell) \wedge t =
 F_{\varphi,1}(s,\bar y) \rightarrow \cT(t) \wedge s \le_{\cT} t,$
\sn
\item[$\bullet$]  $\cT(s) \wedge xRs \wedge \varphi(x,\bar y) \wedge
  \bigwedge\limits_{\ell < \ell g(\bar y)} P(y_\ell) \wedge t
  = F_{\varphi,1}(s,\bar y) \rightarrow xRt,$
\sn
\item[$\bullet$]  $\cT(s) \wedge x'Rs \wedge 
x'RF_{\varphi,1}(s,\bar y) \rightarrow \varphi(x',\bar y)$.
\end{enumerate}

This suffices.  The proof when $\iota = 2$ is similar.

\underline{Clause (j)}:  Similarly but we give details. 

Let $\varphi = \varphi(x,\bar y) \in \Delta$, so the following are Horn:

\begin{enumerate}
\item[$\bullet$]  $\varphi(x_1,\bar y) \wedge P(x_1)
  \wedge \bigwedge\limits_{\ell < \ell g(\bar y)} P(y_\ell) \wedge s = 
F_{\varphi,2}(t,\bar y) \rightarrow s \le_{\cT} t,$

\sn
\item[$\bullet$]  $\varphi(x_1,\bar y) \wedge P(x_1)
  \wedge \bigwedge\limits_{\ell < \ell g(\bar y)} P(y_\ell) \wedge s = 
F_{\varphi,2}(t,\bar y) \rightarrow (\exists x)(xRs \wedge 
\varphi(x,\bar y)),$
\sn
\item[$\bullet$]  $P(x) \wedge \bigwedge\limits_{\ell < \ell g(\bar y)}
P(y_\ell) \wedge s = F_{\varphi,2}(t,\bar y) \wedge z \le_{\cT} t \wedge xRz
  \wedge \varphi(x,\bar y) \rightarrow z \le_{\cT} s$.
\end{enumerate}

\underline{Clause (k)}:  As $\theta = \aleph_0$ this is empty.

This suffices.
\end{PROOF}

\begin{claim}
\label{c10}
Also for $\theta > \aleph_0$ (see \ref{h4}(2)) Claim \ref{c7} 
holds but some of the formulas are in $\bbL_{\theta,\theta}$.
\end{claim}

\begin{PROOF}{\ref{c10}}
\underline{Clause $(k)$}:  When $\theta > \aleph_0$.  

Should be clear
because for each limit ordinal $\delta < \kappa$, we have that the sentences
$\psi_\delta = $ $ (\forall x_0,\dotsc,x_\alpha,\dotsc,x_\delta)(\exists
y)(\forall z) \big( (\bigwedge\limits_{\alpha < \beta < \delta} x_\alpha
\le_{\cT} x_B \le_{\cT} y \le_{\cT} x_\delta) \wedge 
(\bigwedge\limits_{\alpha < \beta < \delta} x_\alpha \le_{\cT}
x_\beta \le_{\cT} z \le_{\cT} y \le_{\cT} x_\delta 
\rightarrow y=z ) 
\big)$ is a Horn sentence
and it expresses ``any $\le_{\cT}$-increasing chain of length $\delta$ has a
$\le$-lub". 
\end{PROOF}

\begin{conclusion}  
\label{c13}
1) Assume

\begin{enumerate}
\item[$(a)$]   $D$ be a filter on $I,$
\sn
\item[$(b)$]  $N$ a model, $\lambda = \|N\| + |\tau_N|,\Delta$ the set
of atomic formulas (in $\bbL(\tau_N)$),
\sn
\item[$(c)$]  $\cT = (\cT,\le_{\cT}) := ({}^{\omega >}\lambda,
\trianglelefteq)^I/D$,
\sn
\item[$(d)$]   $\kappa = \gp^*_{\cT} = \min\{\gt_{\cT},\gp_
\theta (\cT_1) \} $  
  see Definition \ref{h2}(6).
\end{enumerate}

\Then \, the reduced power $N^I/D$ is $(\kappa,1,\Delta)$-saturated.

2) Assume\footnote{Note that   
$ \kappa $ here may be bigger than in part (1).}

\begin{enumerate}
\item[$(a)$]   $D$ is a $\theta$-complete filter on $I,\theta =
  \cf(\theta) > \aleph_0$,
\sn
\item[$(b)$]  $N$ is $(\theta,\Delta)$-saturated, $\Delta$ a set of
atomic formulas,
\sn
\item[$(c)$]   $\cT_1 := ({}^{\theta >}\lambda,\trianglelefteq)^I/D$
\sn
\item[$(d)$]   $\kappa = \min\{\gt_{\cT_1},\gp_
\theta (\cT_1) \}. $ 
\end{enumerate}

\Then \, $N^I/D$ is $( \kappa, \theta, 1, \Delta )$-saturated.

3) We can above replace $N^I/D$ by $N^{\gB}/D$ where $D$ is a filter
on the complete Boolean algebra $\gB$ which has
$(<  \theta )$-distributivity when $\theta > \aleph_0$.


\end{conclusion}

\begin{PROOF}{\ref{c13}}
1) Let $\theta = \aleph_0$ and $\mathbf r_0 = (M_0,\Delta)$ be $\mathbf
r^\theta_{N,\Delta}$ from \ref{h6}, so $\theta_{\mathbf r_0} = \theta$.

By Claim \ref{h6}, $M_0$ is an $\RSP$ hence by Claim \ref{c7} 
also $M = M^I_0/D$ is an $\RSP$.  Now apply the Main Claim \ref{h8}(1).

2) Similarly using \ref{h8}(2).

3) Similarly.
\end{PROOF}

\begin{remark}  
\label{c16}
1) No harm in assuming $\Delta = \{Q(\bar y):Q$ a predicate$\}$.  Note
that allowing bigger $\Delta$ is problematic except in trivial
   cases ($\varphi$ and $\neg \varphi$ are equivalent to Horn
   formulas), see proof of clauses (i),(j) of Definition \ref{h4}.

2) Using \ref{c13}(1) above, if $D$ is an ultra-filter, not surprisingly we get
\cite[Ch.VI,2.6]{Sh:c}, i.e. the theory of dense linear orders is
$\trianglelefteq$-maximal (well, using the translation from dense linear
 orders to trees in \ref{c4}(2) equivalently
\cite[Ch.VI,2.7]{Sh:c}).  The new point here  
is that
\ref{c13} does this also for reduced powers, i.e. for $D$ a filter.

3) So a natural question is: can we replace the strict property by
   $\SOP_2$?  
   We shall show that for 
   reduced power we also have non-peculiar cuts, see
   \S4.

4) Why is the reduced power of a tree not necessarily a tree? 
 Let $ M $ be the tree $({}^{ \omega > }\omega , \triangleleft ) $. Let
$\eta_1 \triangleleft \eta_2 \triangleleft \eta_3 \in {}^{\omega
  >}\omega$ and let $A_1,A_2 \in D^+$ be disjoint and define
$f_\ell:I \rightarrow {}^{\omega >}\omega$ for $\ell=1,2,3$ by:

\begin{enumerate}
\item[$\bullet$]  $f_3(s) = \eta_3$ for $s \in I,$
\sn
\item[$\bullet$]  $f_2(s)$ is $\eta_2$ if $s \in A_2$ and $\eta_0$
  otherwise,
\sn
\item[$\bullet$]  $f_1(s)$ is $\eta_2$ if $s \in A_1$ and $\eta_0$
  otherwise.
\end{enumerate}

\underline{Clearly} if   $N = M^I/D$  then in $ N $   we have:

\begin{enumerate}
\item[$\bullet$]  $f_1/D \triangleleft f_3/D,$
\sn
\item[$\bullet$]  $f_2/D \triangleleft f_3/D,$
\sn
\item[$\bullet$]  $\neg(f_1/D \triangleleft f_2/D),$
\sn
\item[$\bullet$]  $\neg(f_2/D \triangleleft f_1/D),$
\sn
\item[$\bullet$]  $\neg(f_1/D = f_2/D)$.
\end{enumerate}
\end{remark}

\begin{conclusion}
\label{c20}
$N^I/D$ is $(\kappa,1,\Delta_1)$-saturated and $\kappa \ge \theta$ \when \,:

\begin{enumerate}
\item[$(*)$]  $(a) \quad D $ 
is a
  $\theta$-complete filter on $I,$
\sn
\item[${{}}$]  $(b) \quad \Delta \subseteq \{\varphi:
\varphi(x,\bar y) \in \bbL_{\theta,\theta}(\tau_N)$ is atomic (hence
$\in \bbL(\tau_N))\},$
\sn
\item[${{}}$]  $(c) \quad \Delta_1 = c \ell_{< \theta}(\Delta) =$ the
  closure of $\Delta$ under conjunction of $< \theta$ formulas,

\item[${{}}$]  $(d) \quad N$ is $(\theta,\Delta)$-saturated, i.e. 
if $p(x) \subseteq \Delta(N)  = \{\varphi ( x, \bar{ a }) :  
\varphi (x, \bar{ y } )\in \Delta , \bar{ a } \in {}^{ \lg(\bar{ y } )}
M\} $ has cardinality $< \theta$ and is 
finitely satisfiable in $N$ \then \, $p$ is realized in $N,$
\sn
\item[${{}}$]  $(e) \quad \kappa = \min\{\gp_{\cT},\gt_
\theta (\cT) \} $
where $\cT = ({}^{\theta
  >}\lambda,\trianglelefteq)^I/D$ and $\lambda = {}^{\theta >}(\| N \| + 
|\Delta|)$.
\end{enumerate}
\end{conclusion}

\begin{PROOF}{\ref{c20}}
Let $\mathbf r = \mathbf r^\theta_{N,\Delta_1}$ recalling Definition
\ref{h6} and $M_0 = M_{\mathbf r}$.  

Now apply \ref{h8}(2) noting that:

\begin{enumerate}
\item[$(*)_1$]  $N_1 = N^I_0/D$ satisfies: every set of $< \theta$
  formulas from $\Delta(N)$ which is finitely satisfiable in $N_1$ is
  realized in $N_1$.
\end{enumerate}

[Why?  Let $\langle \varphi_\alpha(x,f_{\alpha,0}/D,\dotsc,
f_{\alpha,n(\alpha)-1}/D):\alpha < \alpha_*\rangle$ be finitely 
satisfiable in $N_1$ and $\alpha_* < \theta,\alpha < \alpha_*
\Rightarrow \varphi_\alpha \in \Delta$.  
For every finite $u \subseteq \alpha_*$ we have $N_1 \models
(\exists x)( \bigwedge\limits_{\alpha \in u}
\varphi_\alpha(x,f_{\alpha,0}/D,\ldots,  
f_{\alpha , n(\alpha )-1}/D)\big) $ 
hence the set 

\[
I_u := \{s \in
I:N_1 \models (\exists x) \bigwedge\limits_{\alpha \in u}
\varphi_\alpha(x,f_{\alpha,0}(s),\dotsc,f_{\alpha,n(\alpha)-1}(s))\}
\]

belongs to $D$.  But $D$ is $\theta$-complete, hence $I_* =
\cap\{I_u:u \subseteq \alpha_*$ is finite$\}$ belongs to $D$.  Now for
each $s \in I_*$, the set $p_s :=
\{\varphi_\alpha(x,f_{\alpha,0}(s),\dotsc,f_{\alpha,n(\alpha)-1}(s)):\alpha
< \alpha_*\}$ is finitely satisfiable in $N$, hence is realized by
some $a_s \in N$.  Let $g \in {}^I N$ be such that $s \in I_*
\Rightarrow g(s) = a_s$; clearly $g/D$ realizes $p$, so we are done.]

Similarly

\begin{enumerate}
\item[$(*)_2$]  in $\cT = ({}^{\theta >}\lambda,\trianglelefteq)^I/D$ we have,
\sn
\item[${{}}$]  $(a) \quad$ every increasing sequence of length $<
  \theta$ has an upper bound,
\sn
\item[${{}}$]  $(b) \quad$ any increasing sequence of length $<
 \theta$ with an upper bound has a lub,
\sn
\item[${{}}$]  $(c) \quad$ there is no infinite decreasing sequence so
  $(\kappa_1,\kappa_2) \in \cC ({\cT}) \Rightarrow \kappa_2 = 1$.
\end{enumerate}

[Why?  For clause (a) note that $(\forall x_0,\dotsc,x_\alpha,\ldots)_{\alpha
  < \delta}(\exists y)(\bigwedge\limits_{\alpha <\beta < \delta}
x_\alpha \le_{\cT} x_\beta \rightarrow \bigwedge\limits_{\alpha < \delta}
x_\alpha \le_{\cT} y)$ is a Horn sentence.  For clause (b) see
\ref{c10}, i.e. proof of clause (k) in \ref{c10}.]

\begin{enumerate}
\item[$(*)_3$]  $M_1 = M^I_{\mathbf r}/D$ is a $\theta-\RSP$.
\end{enumerate}

[Why?  See above recalling \ref{c7}, \ref{c10}.]

\begin{enumerate}
\item[$(*)_4$]  if $ \theta \ge   
{\aleph_1} $
then $\mathbf r$ satisfies $(k)^+$ from \ref{h8}(3).
\end{enumerate}

[Why?  Easily as $D$ is a $\aleph_1$-complete ultra-filter.]

So we are done by \ref{h8}(3).
\end{PROOF}

It is natural to wonder
\begin{question}
\label{c24}
Assume $ \lambda \ge \theta = \cf(\theta ) > {\aleph_0} $.

1) Is there a
$\theta$-complete $(\lambda,\theta)$-regular ultra-filter
 $ D $ on
$\lambda$ such that $\lambda < \gt(({}^{\theta>}\theta,
\trianglelefteq)^\lambda/D)$?

2) Similarly for filters.

3) Use $\le_{\cT} = \trianglelefteq$ or $<_{\cT} = \triangleleft$?

4) If $\lambda = \lambda^{< \theta},D$ a fine normal ultra-filter on
   $I=[\lambda]^{< \theta}$, we get $\lambda \le
   \gt({}^{\theta>}\theta,\trianglelefteq)/D$. 
\end{question}

\begin{remark}
\label{c26}
Now \cite[\S5]{Sh:1030} answers \ref{c24}(1) positively for $\theta$
a super-compact cardinal.
\end{remark}

\begin{conclusion}
\label{c28}
Let $\gB$ be a complete Boolean algebra and $D$ a filter on $\gB$.

1) For every model $N$, letting $\lambda = \|N\| + |\tau_N|$, we
have $N^{\mathfrak{B} }/D$ is 
($\mu^+$,atomic)-saturated if $\mu^+ \le
\min\{\gp(({}^{\omega >}\lambda,\trianglelefteq)^{\mathfrak{B} 
}/D)),\gt(({}^{\omega
   >}\lambda,\trianglelefteq)^{\gB}/D)\}$.

2) Assume $\gB$ is 
$(< \theta)$-distributive


(e.g. 
for some  dense  $Y \subseteq \gB^+$, for 
every decreasing sequence\footnote{can weaken the
   demand.} in $\mathfrak{B} $ of elements from $Y$ of length $< \theta$ has a
   positive lower bound),  
and $D$ is a $\theta$-complete filter on $\gB$.  
If $N$ is ($\mu^+$, atomic)-saturated \then \, 
$N^{\mathfrak{B} }/D$ is $\gt(({}^{\theta >}\lambda,
\trianglelefteq)^{\mathfrak{B} }/D)$-atomic saturated. 
\end{conclusion}

\begin{PROOF}{\ref{c28}}
As, e.g. in \ref{c20} above or \ref{c46} below.
\end{PROOF}

\begin{conclusion}\label{c31}
Assume $(T,\varphi(\bar x,\bar y))$ has $\SOP_3$.

Then, recalling \ref{z16}, $T$ is
$\trianglelefteq^{\rp}_\lambda$-maximal for every $\lambda$ and even
$(T,\{\varphi(\bar x,\bar y)\})$ is.
\end{conclusion}

\begin{PROOF}{\ref{c31}}
    By \cite[Ch.~VI.2.7]{Sh:a} $=$ \cite{Sh:c}, but we elaborate.

    Without loss of generality assume $\varphi$ is an atomic formula. 

    Let $(T_{1}, \Delta_{1}) = (T_{1}, \{ \varphi(x, y) \})$ and let $T_{2}$ be a first-order theory, $\Delta_{2}$ a set of atomic formulas in $\bbL(\tau(T_{2})),$ $\vert T_{2} \vert \leq \lambda,$ $\cB = \cP(\lambda),$ $D$ a regular filter on $\lambda,$ or, equivalently, on the Boolean algebra $\cB.$ Finally, let $M_{\ell}$ be a $\lambda^{+}$-saturated model of  $T_{\ell}$ for $\ell = 1, 2.$ We assume $M_{1}^{\lambda} / D$ is $(\lambda^{+}, \Delta_{1})$-saturated and we should prove that also $M_{2}^{\lambda} / D$ is $(\lambda^{+}, \Delta_{2})$-saturated. By Conclusion \ref{c28}, 

    \begin{enumerate}
        \item[$(\ast)_{1}$] It suffice to prove $\mu^{+} \leq \min \{ \mathfrak{p}( ({}^{\omega >} \lambda, \unlhd)^{\gB} / D), \mathfrak{t}( ({}^{\omega >} \lambda, \unlhd)^{\gB} / D) \}.$
    \end{enumerate}

    Next, 

    \begin{enumerate}
        \item[$(\ast)_{2}$] Let $I$ be the following linear order: 
    
        \begin{itemize}
            \item the set of elements is $\{ (\eta, \iota) \colon \eta \in {}^{\omega > } \lambda, \iota \in \{ 1, -1 \}  \},$
    
            \item the linear  order $<_{I}$ is: $(\eta_{1}, \iota_{1}) < (\eta_{1}, \iota_{2})$ iff one of the following occurs: 
    
            \begin{enumerate}
                \item[(a)] $\eta_{1} \lhd \eta_{2}$ and $\iota_{1} = -1,$
    
                \item[(b)] $\eta_{2} \lhd \eta_{1}$ and $\iota_{1} = 1,$
    
                \item[(c)] $\eta_{1} = \eta_{2}$ and $\iota_{1} = -1,$ $\iota_{2} = 1,$
    
                \item[(d)] for some $k < \lg(\eta_{1}),$ $\lg(\eta_{1}),$ we have $\eta_{1} \, {\rest} \, k= \eta_{2} \, {\rest} \, k$ and $\eta_{1}(k) < \eta_{2}(k),$ 
            \end{enumerate}
        \end{itemize}
    \end{enumerate}
    
    Next, we can find $a_{(\eta, \iota)} \in M_{1}$ for $(\eta, \iota) \in I$ such that $(\eta_{1}, \iota_{1}) <_{I} (\eta_{2}, \iota_{2}) \Rightarrow M_{1} \models \varphi[a_{(\eta_{1}, \iota_{1})}, a_{(\eta_{2}, \iota_{2})}].$

    Let $N = (\vert M_{1} \vert, P^{N}, R^{N}),$ where: $$P^{N} = \{ a_{(\eta, \iota)} \colon (\eta, \iota) \in I \} \text{ and } R^{N} = \{ 
    (a_{(\eta_{1}, \iota_{1})}, a_{\eta_{2}, \iota_{2}}) \colon (\eta_{1}, \iota_{1}) <_{I} (\eta_{2}, \iota_{2})\}.$$ 

    Let $M_{1}^{+} = (M_{1}, P^{N}, R^{N}).$ Now, 

    \begin{enumerate}
        \item[$(\ast)_{3}$] in $(M_{1}^{+})^{\lambda / D}$ every set of formulas $\subseteq \{ 
        \varphi(x, a), \varphi(a, x) \colon a \in (M_{1}^{+})^{\lambda} / D \}$ of cardinality $\leq \lambda,$ which is finitely satisfiable is realized.  
    \end{enumerate}

    [Why? Think of $M_{1}^{\lambda} / D$.]

    \begin{enumerate}
        \item[$(\ast)_{4}$] In $N^{\lambda} / D$ every set  of formulas $\subseteq \{ x R a, aRx \colon a \in N^{\lambda} / D \}$ of cardinality $\leq \lambda,$ which is finitely satisfiable is realized.
    \end{enumerate}

    So easily we are done. 
\end{PROOF}

\centerline {$* \qquad * \qquad *$}

On the connection to Peano arithmetic and to Pabion \cite{Pab82}, see
Malliaris-Shelah \cite{Sh:1051}.  We repeat some results of
\cite{Sh:1069} in the present context; but first recalling:

\begin{definition} \label{c39}
1) \PA, Peano arithmetic,  is the f.o. theory consisting of: 
\begin{enumerate}
\item[(a)] the obvious axioms on $ 0,1,x < y , x + y , x y $
\item[(b)] all the cases of the induction scheme, i.e.  for every
f.o. $ \varphi $:

``if $ \{x: \varphi (x, \bar{ y }  )  \} $ is not empty
then it has a first member", 
\end{enumerate}

2) \BPA, the bounded  Peano arithmetic, is defined similarly,
but in clause (b), the formulas $ \varphi $ is bounded, i.e.
all the quantifications inside it are of the form
$ (\forall x < y )$ or $ (\exists x < y  ) $.
\end{definition}

\begin{definition}
\label{c42}

%
1) $N \models \BPA$ is boundedly $\kappa$-saturated up to
$(c_1,c_2)$ where $c_1,c_2\in N$ \when \,: if $p(x) 
\cup \{x < c_1\}$ is a type in $N$ (= finitely satisfiable) of
  cardinality $< \kappa$ consisting of bounded formulas but with
parameters $\le c_2$, \then \, $p(x) \cup \{x < c_1\}$ is realized in
$N$.

2) If above $c_1 = c = c_2$ we may write $c$ instead of $(c_1,c_2)$.  We say $N$ is strongly
boundedly $\kappa$-saturated up to $c$ \when \, it holds for $(c,c_2),c_2
= \infty$, i.e. we do not bound the parameters.

3) Omitting ``up to $c$" in part (3) means for every $c \in N$.
\end{definition}

\begin{conclusion}
\label{c46}
Assume $N$ is a model of $\BPA$.

1) Assume $a_* \in N$ is non-standard and the power in the $N$-sense 
$c^{a_*}$ exists for every $c \in N$.

For any uncountable cardinal $\kappa$ the following conditions are
equivalent:

\begin{enumerate}
\item[$(a)$]  $N$ is boundedly $\kappa$-saturated up to $c$ for any $c
  \in N,$
\sn
\item[$(b)$]  if $(C_1,C_2)$ is a cut of $N$ of
cofinality $(\kappa_1,\kappa_2)$ and $\kappa_1,\kappa_2$ are infinite (so
$C_1,C_2 \ne 0$) \then \, $\kappa_1 + \kappa_2 \ge \kappa,$
\sn
\item[(c)]  like clause (b) but $ \kappa_1=\kappa_2$,
 that is restricting ourselves to symmetric cuts.  
\end{enumerate}

2) We can weaken the assumption of part (1) by fixing $c$, as well as
$N,a_*$.  That is, assume $N \models ``n < a_*$ and $c_n =
c^{(a_*)^n}$ exist" for every standard $n$ from $N$.  For 
every uncountable cardinal $\kappa$ the following are
equivalent:

\begin{enumerate}
\item[$(a)'$]  $N$ is boundedly $\kappa$-saturated up to $c_n$ for
each $n$
\sn
\item[$(b)'$]  if $(C_1,C_2)$ is a cut of $N$ of cofinality
$(\kappa_1,\kappa_2)$ with $\kappa_1,\kappa_1$ infinite such that
$c_n \in C_2$ for some $n$ \then \, $\kappa_1 +
\kappa_2 \ge \kappa$
\sn
\item[(c)']  like clause (b)' but $ \kappa_1=\kappa_2$.
\end{enumerate}

3) Moreover we can add in part (2):

\begin{enumerate}
\item[$(c)$]  $N$ is strongly boundedly $\kappa$-saturated up to $c$.
\end{enumerate}
\end{conclusion}

\begin{PROOF}{\ref{c46}}
1) By (2).

2) \underline{$(a)' \Rightarrow (b)'$}: 

Trivial.
\smallskip

\underline{$(b)' \Rightarrow (a)'$}: 

\Wilog \, $c$ is not standard (in $N$) and $n=0$.
Let $N^+ = (N,c,a_*)$ and $\tau^+ = \tau(N^+) = 
\tau(N) \cup \{c,a_*\}$ and $\Delta = \{\varphi(x,\bar y) \wedge 
x < c \wedge \bigwedge\limits_{\ell} y_\ell < c:\varphi(x,\bar y) \in
\bbL(\tau_N)$ is a bounded formula$\}$.  We define $\mathbf r$ naturally -
the tree of sequences of length $< a_*$ of members of $\Delta(N_{\le c})$
possibly non-standard but of length $< a_*$.  Now apply \ref{h8}.

\underline{$(b)' \Rightarrow (c)'$}: 

Obvious.

\underline{$(c)' \Rightarrow (b)'$}: 

By \cite{Sh:998}.

3) We just repeat the proof of \ref{h8} or see \ref{c59} below.
\end{PROOF}

\begin{question}
\label{c48}
Is $a_*$ necessary in \ref{c46}(1)?  We conjecture that yes.
\end{question}

A partial answer:
\begin{fact}
\label{c52}
If $N$ is a model of $\PA$, \then \, $N$ is $\kappa$-saturated \Iff \,
$\cf(|N|,<^N) \ge \kappa$ and $N$ is boundedly $\kappa$-saturated. 
\end{fact}



%

\begin{claim}
\label{c59}
If (A) then (B) where:

\begin{enumerate}
\item[$(A)$]  $(a) \quad \mathbf r_\alpha$ is an $\RSP$ for $\alpha <
  \delta,$
\sn
\item[${{}}$]  $(b) \quad \Delta_{\mathbf r_\alpha} = \Delta$ is a set
  of quantifier-free formulas,
\sn
\item[${{}}$]  $(c) \quad \cT_{\mathbf r_\alpha} = \cT_{\mathbf r_0}$ and
  $N_{\mathbf r_\alpha}$ is increasing with $\alpha,$
\sn
\item[${{}}$]  $(d) \quad Q \in \tau(N_{\mathbf r_\alpha})$ and
 $Q^{N_{\mathbf r_\alpha}} = Q^{N_{r_0}},$
\sn
\item[${{}}$]  $(e) \quad$ if $\varphi(x,\bar y) \in 
\Delta_{\mathbf r_\alpha}$ and $\bar b \in {}^{\ell g(\bar y)}
(N_{\mathbf r_\alpha})$ \then \, $\varphi(N_{\mathbf r_\alpha},\bar b)
\subseteq Q^{N_{\mathbf r_\alpha}},$
\sn
\item[${{}}$]  $(f) \quad \kappa = \min\{\gp_\alpha(\cT_{\mathbf
    r_0}),\gt(\cT_{\mathbf r_0})\}.$ 
\sn
\item[$(B)$]  the model $\cup\{N_{\mathbf r_\alpha}:\alpha < \delta\}$
  is $(\kappa,1,\Delta)$-saturated.
\end{enumerate}
\end{claim}

\begin{PROOF}{\ref{c59}}  
As in \ref{h8}.
\end{PROOF}

\newpage

\section {Criterion for atomic saturation of reduced powers} \label{3}\

Malliaris-Shelah \cite{Sh:998} have  dealt with such problems for
ultra-filters (on sets).  The main case here is $\theta = \aleph_0$.
\begin{definition}
\label{g3}
Assume $D$ is a filter on the complete Boolean algebra $\gB,T$ an
$\bbL_{\theta,\theta}(\tau_T)$-theory, $\Delta \subseteq \bbL(\tau_T)$
and $\mu \ge |\Delta|$.  We say
$D$ is a $(\mu,\theta,\varepsilon!,\Delta,T)$-moral filter on $\gB$
(writing $\varepsilon$ instead $\varepsilon$! 
means for every $\varepsilon' < 1 +\varepsilon$;
if $\gB = \cP(\lambda)$ we may say good instead of moral): \when \, for every
$D-(\mu,\theta,\varepsilon!,\Delta,T)$-problem there is a
$D-(\mu,\theta) $ -
solution \underline{where}:

\begin{enumerate}
\item[$(a)$]  $\bar{\mathbf a}$ is a 
$D-(\mu,\theta,\varepsilon !,\Delta,T)$-(moral)-problem \when \,:
\sn
\begin{enumerate}
\item[$(\alpha)$]  $\bar{\mathbf a} = \langle \mathbf a_u:u \in [\mu]^{<
  \theta}\rangle,$
\sn
\item[$(\beta)$]  $\mathbf a_u \in D$ (hence $\in \gB^+$),
\sn
\item[$(\gamma)$]  $\bar{\mathbf a}$ is $\subseteq$-decreasing, that is
  $u \subseteq v \in [\mu]^{< \theta} \Rightarrow \mathbf a_v \le \mathbf
  a_u$ and $\mathbf a_\emptyset = 1_{\gB},$
\sn
\item[$(\delta)$]  for some sequence $\langle 
\varphi_\alpha(\bar x_{[\varepsilon]},\bar y_\alpha):\alpha 
< \mu\rangle$ of formulas
from $\Delta$ for every $\mathbf a \in \gB^+ $ 
and  
$u \subseteq \mu$ of cardinality $< \theta$ we can find $M
  \models T$ and $\bar b_\alpha \in {}^{\ell g(\bar y_\alpha)}M$ 
for $\alpha \in u$ such that: 
\sn
\item[${{}}$]  $\bullet \quad$ for every $v \subseteq u$ we have:
 
 $\mathbf a \le \mathbf a_v \Rightarrow M \models 
``(\exists \bar x_{[\varepsilon]}) \, 
\bigwedge\limits_{\alpha \in v} \varphi_\alpha(\bar
x_{[\varepsilon]},\bar b_\alpha)"$

\hskip35pt  and $\mathbf a \le 1 - \mathbf a_v 
\Rightarrow M \models ``\neg(\exists \bar x_{[\varepsilon]})
  \bigwedge\limits_{\alpha \in v} \varphi_\alpha(\bar
  x_{[\varepsilon]},\bar b_\alpha)"$
\end{enumerate}
\sn
\item[$(b)$]  $\bar{\mathbf b}$ is a $D-(\mu,\theta)$-(moral)-solution 
of the $D-(\mu,\theta,\varepsilon!  
,\Delta,T)$-(moral)-problem 
$\bar{\mathbf a}$ \when \,:  
\sn
\begin{enumerate}
\item[$(\alpha)$]  $\bar{\mathbf b} = \langle \mathbf b_u:u \in 
[\mu]^{< \theta}\rangle,$
\sn
\item[$(\beta)$]  $\mathbf b_u \in D$ and $\mathbf b_\emptyset = 1_{\gB},$
\sn
\item[$(\gamma)$]  $\mathbf b_u \le \mathbf a_u,$
\sn
\item[$(\delta)$]  $\bar{\mathbf b}$ is multiplicative,
  i.e. $\bar{\mathbf b}_u = \cap\{\mathbf b_{\{\alpha\}}:\alpha \in u\}$ 
  and $ \mathbf{b} _ \emptyset = 1_{\mathfrak{B} }$.
\end{enumerate}
\end{enumerate}
\end{definition}

\begin{remark}
\label{g5}
1) The $\theta$ here means ``a type is $(< \theta)$-satisfiable".

2) The use of ``$\varepsilon!$" is to 
conform with Definition \ref{z6}.
\end{remark}

Recall (from \ref{z6}):

\begin{definition}
\label{g8}
1) 
Let $\tau$ be a vocabulary and $\Delta \subseteq \{\varphi \in
\bbL(\tau):\varphi = \varphi(\bar x,\bar y)\}$ but $\varphi(\bar x,\bar y) \in \Delta$ means we can add to
$\bar x$ dummy variables.  Let $\lambda > \theta$ (dull otherwise).

A $\tau$-model $M$ is $  
(\lambda,\theta,\varepsilon !,
\Delta )$
-saturated
\when \,: if $p \subseteq \{\varphi(\bar x_{[\varepsilon]},\bar
a):\varphi(\bar x_{[\varepsilon]},\bar y) \in \Delta,\bar a \in {}^{\ell
  g(\bar y)}M\}$ has cardinality $< \lambda$ 
  and 
is $(< \theta)$ satisfiable in $M$ \then \, $p$ is realized in $M$. 
\end{definition}

\begin{claim}
\label{g9}
1) For a $(\mu,\theta)$-regular $\theta$-complete ultra-filter 
$D$ on a set $I$ and $\theta$-saturated or just
$(\theta,\aleph_0,\varepsilon !,\Delta)$-saturated model $M$, a
cardinal $\mu$ and $\Delta =
\bbL_{\theta,\theta}(\tau_M)$, the following conditions
are equivalent:

\begin{enumerate}
\item[$(a)$]  $D$ is $(\mu,\theta,\varepsilon !,\Delta,T)$-moral ultra-filter on the
  Boolean algebra $\cP(I),$
\sn
\item[$(b)$]  if $M \in \Mod_T$ then $M^I/D$ is 
$(\mu,\theta,\varepsilon !,\Delta)$-saturated.
\end{enumerate}

2) Similarly for $D$ a ultra-filter on a $(< \theta  ) $-distributive (see \ref{z9}(8))  
complete Boolean algebra $\gB$.
\end{claim}

\begin{PROOF}{\ref{g9}}
    Similar to \ref{g10}, it actually follows from it because as $D$ is an
    ultra-filter, we can start with $M \models T$, expand it to $M^+$ by
    adding a predicate to any definable relation and applying \ref{g10} to
    $T^+ = \Th(M^+)$.
\end{PROOF}

\begin{claim}
\label{g10}
1) If $(A)$ then $(B) \Leftrightarrow (C)$ where:

\begin{enumerate}
\item[$(A)$]  $(a) \quad \gB = \cP(I),$
\sn
\item[${{}}$]  $(b) \quad D$ is a $\theta$-complete
$(\mu,\theta)$-regular filter on $\gB,$
\sn
\item[${{}}$]  $(c) \quad $ 
$\theta > \varepsilon$ or just $ \mu ^+ > \varepsilon,$  
\sn
\item[${{}}$]  $(d) \quad T$ is an $\bbL_{\theta,\theta}(\tau)$-theory,
\sn
\item[${{}}$]  $(e) \quad \Delta$ is a set of conjunctions of $<
\theta$ atomic formulas from $\bbL_{\theta,\theta}(\tau).$
\sn
\item[$(B)$]  $D$ is a 
$(\mu,\theta,\varepsilon !,\Delta,T)$-moral filter on $\gB.$ 
\sn
\item[$(C)$]   if $M_s$ is a model of $T$ for $s \in I$
\then \, $\prod\limits_{s \in I} M_s/D$ is
$(\mu^+,\theta,\varepsilon !,\Delta)$-saturated.
\end{enumerate}

2) If $(A)'$ then $(B)' \Leftrightarrow (C)'$ where:

\begin{enumerate}
\item[$(A)'$]  $(a) \quad \gB$ is a $(< \theta)$-distributive (see
  \ref{z9}(8)) complete Boolean algebra,
\sn
\item[${{}}$]  $(b)-(e) \quad$ as above (on regularity see Definition,
  \ref{z10}) 
\sn
\item[${{}}$]  $(d)^+ \quad T$ is a complete
  $\bbL_{\theta,\theta}(\tau)$-theory.
\sn
\item[$(B)'$]  as $(B)$ above.
\sn
\item[$(C)'$]  $(a) \quad$ if $M$ is a model of $T$ \then \, $M^{\gB}/D$ is
$(\mu^+,\theta,\varepsilon !,\Delta)$-saturated,
\sn
\item[${{}}$]  $(b) \quad$ if $\cI$ is a maximal anti-chain of $\gB$
  and $\bar M = \langle M_b:b \in \cI\rangle$ is a 

\hskip25pt  sequence of $\tau$-models then $\bar M^{\gB}/D$ is $(\mu^+,\theta,\varepsilon
  !,\Delta)$-saturated. 
\end{enumerate}
\end{claim}

\begin{PROOF}{\ref{g10}}
1) \underline{Proving $(B) \Rightarrow (C)$}:
Let $N = \prod\limits_{s \in I} M_s/D$ let 
$\bar x = \bar x_{[\varepsilon]},\varphi_\alpha 
= \varphi_\alpha(\bar x,\bar y_\alpha)$ and assume that
$p(\bar x) = \{\varphi_\alpha(\bar x,\bar b_\alpha):\alpha < \alpha_*\}$ is
$(< \theta)$-satisfiable in $N$ and $|\alpha_*| \le \mu$, so
\wilog \, $\alpha_* = \mu$; \wilog \, let $\varphi_\alpha = \varphi_\alpha(\bar
x,\bar y_{[\xi_\alpha]})$ so $\bar b_\alpha \in
{}^{\xi_\alpha}(\prod\limits_{s \in I} M_s)$.

Let $\bar b_\alpha = \langle f_{\alpha,\xi}/D:\xi < 
\xi_\alpha\rangle$ where $f_{\alpha,\xi} \in \prod\limits_{s \in I} M_s$ and
for $s \in I$ let $\bar b_{\alpha,s} = \langle f_{\alpha,\xi}(s):\xi < 
\xi_\alpha\rangle$; now for $u \in [\mu]^{< \theta}$ we let

\begin{enumerate}
\item[$(*)_0$]  $\mathbf a_u := \{s \in I:M_s \models (\exists \bar x)
\bigwedge\limits_{\alpha \in u} \varphi(\bar x,\bar b_{\alpha,s})\}$.
\end{enumerate}

Now

\begin{enumerate}
\item[$(*)_1$]  $\bar{\mathbf a} = \langle \mathbf a_u:u \in [\mu]^{<
  \theta}\rangle$ is a $D-(\mu,\theta,\varepsilon !,\Delta,T)$-problem.
\end{enumerate}

[Why?  We should check Definition \ref{g3}, clause (a): now 
$(a)(\alpha)$ is trivial;
also 
$\mathbf a_u \subseteq I$ holds by the choice
of $\mathbf a_\alpha$. 
Toward clause $(a)(\beta)$ fix a set $u \in
[\mu]^{< \theta}$;  some $\bar c \in {}^\varepsilon N$ realizes the type
$p_u(\bar x_{[\varepsilon]}) = \{\varphi_\alpha(\bar x,\bar
b_\alpha):\alpha \in u\}$ in $N$
because $ p(\bar{ x } )$ is $ (< \theta ) $-satisfiable in $ N $, 
see Definition \ref{g8}, so let
$\bar c = \langle g_\zeta/D:\zeta < \varepsilon\rangle$ for some
$g_\zeta \in \prod\limits_{s \in I} M_s$ for $\zeta < \varepsilon$ and
let $\bar c_s = \langle g_\zeta(s):\zeta < \varepsilon\rangle \in
{}^\varepsilon(M_s)$.  So $\mathbf a'_{\{\alpha\}} = \{s \in I:M \models
\varphi_\alpha[\bar c_s,\bar b_s]\}$ belong to $D$ because $N \models
\varphi_\alpha[\bar c,\bar b_\alpha]$ by the definition of $N$ if
$\varphi_\alpha$ is atomic, but recalling $D$ is $\theta$-complete
also for our $\varphi_\alpha$, remembering clause (A)(e) 
of \ref{g10}(1).
As $D$ is $\theta$-complete clearly, 
$\mathbf a'_u = \cap\{\mathbf a'_{\{\alpha\}}:\alpha \in u\}$
belongs to $D$ and by our choices, $\mathbf a'_u \le_{\gB} \mathbf a_u$,
hence $\mathbf a_u \in D$ so subclause $(a)(\beta)$ 
of Def \ref{g3}  
holds indeed.

By the choice of $\mathbf a_u,\bar{ \mathbf{a}}$ is $\subseteq$-decreasing 
so subclause $(a)(\gamma)$  of Def \ref{g3} holds.

Lastly, subclause $(a)(\delta)$ of Def \ref{g3} holds by the definition of $\mathbf
a_u$'s recalling $p(\bar x)$ is $(< \theta)$-satisfiable 
(and $ \emptyset \notin D)$.]

\begin{enumerate}
\item[$(*)_2$]  there is $\bar{\mathbf b}$, a $D-(\mu,\theta)$-solution
of $\mathbf a$ in $\gB$.
\end{enumerate}

[Why?  Because we are presently assuming clause (B) of \ref{g10} which
says that $D$ is $(\mu,\theta,\varepsilon!,\Delta,T)$-good, 
see Definition \ref{g3}.]

\begin{enumerate}
\item[$(*)_3$]  \wilog \, $s \in I \Rightarrow \{\alpha < \mu:s \in
  \mathbf b_{\{\alpha\}}\}$ has cardinality $< \theta$.
\end{enumerate}

[Why?  As $D$ is $(\mu,\theta)$-regular.]

Next for $s \in I$ let $u_s = \{\alpha < \mu:s \in \mathbf
b_{\{\alpha\}}\}$ but $\bar{\mathbf b}$ is multiplicative (see
 \ref{g3}$(b)(\delta)$) so $\mathbf b_{u_s} = \cap\{\mathbf
 b_{\{\alpha\}}:\alpha \in u_s\} = \cap\{\mathbf b_\alpha$: the ordinal
 $\alpha$ satisfies $s \in \mathbf b_{\{\alpha\}}\}$ hence 
$s \in \mathbf b_{u_s}$ hence (see
  \ref{g3}(b) recalling that $|u_\alpha| < \theta$ by $(*)_2$) 
we have $s \in \mathbf a_{u_s}$ hence (by the choice of $\mathbf
  a_{u_s}$) there is $\bar a_s \in {}^\varepsilon(M_s)$ realizing
$\{\varphi(\bar x_{[\varepsilon]},\langle
  f_{\alpha,\varepsilon_\zeta}(s):\zeta < \varepsilon\rangle):
\alpha \in u_s\}$.

Let $\bar a_s = \langle a_{s,\zeta}:\zeta < \varepsilon\rangle$.  
Now for $\zeta < \varepsilon = \ell g(\bar x)$ let $g_\zeta \in
\prod\limits_{s \in I} M_s$ be defined by $g_\zeta(s) =
a_{s,\zeta} \in M_s$ and let $\bar a = \langle
g_\zeta/D:\zeta < \varepsilon \rangle$ noting $g_\zeta/D \in
\prod\limits_{s \in I} M_s/D = N$.  
Hence for every $\alpha < \mu,\{s \in I:M_s \models
\varphi_\alpha(\langle g_\zeta(s):\zeta < \varepsilon \rangle,
\bar b_{\alpha,s})\} \supseteq \mathbf b_{\{\alpha\}} \in D$ so $N
\models \varphi[\bar a,\bar b_\alpha]$.

Hence $\bar a$ realizes $p(\bar x)$ in $N$ as promised.

\underline{Proving $(C) \Rightarrow (B)$}:

To prove clause (B), 
let $\bar{\mathbf a}$ be a $D-(\mu,\theta,\varepsilon
!,\Delta,T)$-problem and let $\bar\varphi = \langle
\varphi_\alpha(\bar x_{[\varepsilon]},\bar y_\alpha):\alpha <
\mu\rangle$ be a sequence of formulas from $\Delta$ as in clause
$(a)(\delta)$ of Definition \ref{g3}.

As $ D $ is $ (\lambda, \theta  ) $-regular, we can choose 
$\bar w = \langle w_s:s \in I \rangle$ 
a sequence of subsets
of $\mu$ each of cardinality $< \theta$ such that $\alpha < \mu
\Rightarrow \{s \in I:\alpha \in w_s\} \in D$.  For $u \in 
[\mu]^{< \theta}$ let $\mathbf c_u = \{s \in I:u \subseteq w_s\}$,
so clearly $\mathbf c_u \in D$ and $\langle \mathbf c_u:u \in [\lambda]^{<
  \theta}\rangle$ is multiplicative.

For each $s \in I$ applying Definition \ref{g3}$(a)(\delta)$ to $\mathbf
a = \{s\}$ and $u = w_s$ we can find a model $M_s$ of $T$ and $\bar
b_{s,\alpha} \in {}^{\ell g(\bar y_\ell)}(M_s)$ for $\alpha \in w_s$
satisying $\bullet$ there.

Now choose $\bar b_{s,\alpha}$ also for $s \in I,\alpha \in \mu
\backslash w_s$, as any sequence of members of $M_s$ of length
$\ell g(\bar y_\alpha)$.  Now for every $\alpha < \mu$ and $j <
\ell g(\bar y_\alpha)$ we define $g_{\alpha,j} \in \prod\limits_{s \in
  I} M_s$ by $g_{\alpha,j}(s) = (\bar b_{s,\alpha})_j$.  

Hence $g_{\alpha,\zeta}/D \in \prod\limits_{s \in I} M_s/D = N$ and
$\bar b_\alpha = \langle g_{\alpha,\zeta}/D:\zeta < \ell g(\bar
y_\alpha)\rangle \in {}^{\ell g(\bar y_\alpha)}N$ and consider the set
$p = \{\varphi_\alpha(\bar x,\bar b_\alpha):\alpha < \mu\}$.  Is
$p$ a $(< \theta)$-satisfiable type in $N$?  
We shall prove that Yes, 
so let  
$u \in
[\mu]^{< \theta}$, then recall
$\mathbf c_u = \{s \in I: u \subseteq W_s \}  \in D$ and 
$s \in \mathbf{c} _u \cap \mathbf{a} _ u  
\Rightarrow \{\varphi_\alpha(\bar x_{[\varepsilon]},\bar
b_{s,\alpha}):\alpha \in u\}$ is realized in $M_s$,
[why? by the choice of 
$ \langle \bar{ b_{s, \alpha }}:  \alpha \in w_s  \rangle $.]


So let the type  $ \{\varphi _ \alpha (\bar{ x } _{[\varepsilon ]},
\bar{ b } _{s, \alpha }):  \alpha \in w_ s  \} $ 
be realized
$\bar a_s =
\langle a_{s,\zeta}:\zeta < \varepsilon\rangle$; for $s \in I $ 
and let $f_{\alpha ,\zeta} \in \prod\limits_{s \in I} M_s$ be
$f_{\alpha,\zeta}(s) = a_{s,j}$. 
Easily $\langle f_{\alpha,j}/D:\zeta
< \varepsilon\rangle$ realizes 
$\{\varphi _ \alpha (\bar{ x } _{[\varepsilon ]}): 
\alpha \in u \} $ 
because $\mathbf a_u \cap \mathbf c_u
\in D$. Hence 
 $ p(\bar{ x } _{[\varepsilon ]})$ is $ (< \theta )$-satisfiable 
 indeed.

Next, we apply clause (C) we are assuming hence
$p(\bar
x_{[\varepsilon]})$ is realized in $N$.
So let $\bar a = \langle a_\zeta:\zeta < 
\varepsilon\rangle \in {}^\varepsilon N$
realize $p$ and let $a_\zeta = h_\zeta/D$ where $h_\zeta \in
\prod\limits_{s \in I} M_s$ and lastly let: 
\[
\mathbf b_u = \{s \in I:M_s \models \varphi_\alpha[\langle
h_\zeta(s):\zeta < \varepsilon\rangle,\bar b_{s,\alpha}] \text{ for
  every } \alpha \in u \text{ and }s \in \mathbf c_u\}.
\]

Now check that $\langle \mathbf b_u:u \in [\lambda]^{< \theta}\rangle$
is as required, recalling $\langle \mathbf c_u:u \in [\lambda]^{<
  \theta}\rangle$ is multiplicative.  So the desired conclusion of
\ref{g3}(B) holds indeed so we are done proving $(C) \Rightarrow (B)$.

2) Similarly; e.g. for clause (a) let $p(\bar x)$ be as there but

\begin{enumerate}
\item[$\bullet$]  $f_{\alpha,\xi} \in M^{\gB}$ is supported by the
maximal anti-chain $\langle \mathbf c_{\alpha,\xi,i}:i <
i(\alpha,\xi)\rangle$
\sn
\item[$(*)_0$]  $\mathbf a_u = \sup\{\mathbf c$: we have $\alpha \in u \wedge \xi <
\xi_\alpha \Rightarrow (\exists \mathbf d)(\mathbf d \in 
\dom(f_{\alpha,\xi}) \wedge \mathbf c \le \mathbf d)$ 
and $M \models (\exists \bar x_{[\varepsilon]}) 
\bigwedge\limits_{\alpha \in u} \varphi(\bar x_{[\varepsilon]},
\langle f_{\alpha,\xi}(\mathbf c):\xi < \xi_\alpha\rangle)\}$
\sn
\item[$(*)_1$]  $\bar{\mathbf a} = \langle \mathbf a_u:u \in 
[\mu]^{< \theta}\rangle$ is a $D-(\mu,\theta,\varepsilon !,\Delta,T)$-problem.
\end{enumerate}

[Why?  As there.]

\begin{enumerate}
\item[$(*)_2$]  let $\bar{\mathbf b}$ be a 
$D-(\mu,\theta)$
-solution.
\end{enumerate}

[Why does $ \bar{ \mathbf{b} } $
exist?  By $(B)'$ recalling Definition \ref{g3}.]

Also the rest is as above.
\end{PROOF}

\begin{remark}
\label{g22}
If $\cS \subseteq [\mu]^{< \theta}$ is cofinal, $u \in [\mu]^{<
  \theta} \Rightarrow |\cP(u) \cap \cS| < \theta_1$ we
  may consistently replace $[\mu]^{< \theta}$ by $\cS$ and
  $2^{\theta_1}$ by $\theta_1$.
\end{remark}

\begin{definition}
\label{g13}
1) A filter $D$ on a complete Boolean algebra $\gB$ is
$(\mu,\theta)$-excellent \when \,: if $\bar{\mathbf a} = \langle \mathbf
   a_u:u \in [\mu]^{< \theta}\rangle$ is a sequence of members of
   $\gB$, (yes! not necessarily from $D$) \then \, we can find 
$\bar{\mathbf b}$ which is a multiplicative
   refinement of $\bar{\mathbf a}$ for $D$, meaning:

\begin{enumerate}
\item[$(a)$]  $\bar{\mathbf b} = \langle \mathbf b_u:u \in [\mu]^{<
\theta}\rangle,$
\sn
\item[$(b)$]  $\mathbf b_u \le \mathbf a_u$ and $\mathbf b_u = \mathbf a_u \mod D,$
\sn
\item[$(c)$]  if $\mathbf a_{u_1} \cap \mathbf a_{u_2} = \mathbf a_{u_1 \cap
u_2} \mod D$ then $\mathbf b_{u_1} \cap \mathbf b_{u_2} = \mathbf b_{u_1
\cap u_2}$.
\end{enumerate}

2) For a Boolean algebra $\gB$ and filter $D$ on $\gB$ we say $\bar{\mathbf
  a}$ is a $D-(\mu,\theta)$-problem
  (or a $ D-(\mu, \theta )$-moral problem)
  when clauses
$(a)(\alpha),(\beta),(\gamma)$ of Definition \ref{g3} holds.

3) A filter $ D $ on a complete Boolean algebra
 $ \mathfrak{B} $ is $ (\mu, \theta  ) $-good 
 \when \, every $ D-( \mu, \theta ) $-problem
 has a $ D-(\mu, \theta )$-solution
\end{definition}

\begin{claim}
\label{g15}
1) Assuming $ (*) $ below, the filter $D$ on $I$ (i.e. on the Boolean algebra $\cP(I)$) is
$(\mu,\theta,\varepsilon !,\Delta,T)$-moral 
\Iff \, 
the filter $D_1$ on 
$\gB_1$ is
$(\mu,\theta,\varepsilon !,\Delta,T)$-moral where:

\begin{enumerate}
\item[(*)] 
\begin{enumerate}  
\item[$(a)$] $ \mathfrak{B} _1 $ is a  complete Boolean algebra,
\sn
\item[$(b)$] $\mathbf j$ is a homomorphism from $\cP(I)$ onto
  $\gB_1$,
\sn
\item[$(c)$]  $D_0 = \{A \subseteq I:\mathbf j(A) = 1_{\gB_1}\}$
  is a $(\mu,\theta)$-excellent filter on $I$, 
\sn
\item[$(d)$] 
$ D_1 $ is a filter on $ \mathfrak{B} _1,$
\sn 
\item[$(e)$] $D= \{A \subseteq I:\mathbf j(A) \in D_1\}$.
is a filter on $ I $
\end{enumerate}
\end{enumerate} 

2) We can replace $\cP(I)$ by a complete Boolean algebra $\gB_2$.
\end{claim}

\begin{PROOF}{\ref{g15}}

\underline{The ``if" direction}:

We assume $D_1$ is $(\mu,\theta,\varepsilon !,\Delta,T)$-moral and
should prove it for $D$.  So let $\bar A = \langle A_u:u \in
[\mu]^{< \theta}\rangle$ be a
$D-(\mu,\theta,\varepsilon !,\Delta,T)$-problem and we should find a
$D-(\mu,\theta)$-solution $\bar B$ of it.

Clearly $\mathbf a_u := \mathbf j(A_u) \in \gB^+$ and $\bar{\mathbf a} =
\langle \mathbf a_u:u \in [\mu]^{< \theta}\rangle = \langle \mathbf
j(A_u):u \in [\mu]^{< \theta}\rangle$ is a
$D_1-(\mu,\theta,\varepsilon!,\Delta,T)$-problem.

Hence by our present assumption ($D_1$ is
$(\mu,\theta,\varepsilon !,\Delta,T)$-moral) there is a
$D_1-(\mu,\theta)$-solution $\bar{\mathbf b}$ of $\bar{\mathbf a}$, let
$\bar{\mathbf b} = \langle \mathbf b_u:u \in [\mu]^{< \theta}\rangle$ so 
in particular 
$u \in [\mu]^{< \theta} \Rightarrow \mathbf b_u \in D_1$.  For $u \in
[\mu]^{< \theta}$ choose $B^1_u \subseteq I $ 
such that $\mathbf
j(B^1_u) = \mathbf b_u$, possible because $\mathbf j$ is a homomorphism
from $\cP(I)$ onto $\gB_1$.  So $\bar B^1 = \langle B^1_u:u \in
[\mu]^{< \theta}\rangle$ is a multiplicative modulo $D_0$, i.e. $\langle
B^1_u/D_0:u \in [\mu]^{< \theta}\rangle$ is a multiplicative sequence of
members of $\cP(I)/D_0$.

Let $B^2_u = B^1_u \cap A_u$, let  

\begin{enumerate}
\item[$\bullet$]  $B^1_u \subseteq A_u \mod D_0$.
\end{enumerate}

[Note that we have written $ B^1_u$ and not $ B^2_u$. 
So 
why this statement holds?  As $\mathbf j(B^1_u) =
\mathbf b_u \le \mathbf a_u = \mathbf j(A_u)$.]

\begin{enumerate}
\item[$\bullet$]  $B^2_u \subseteq B^1_u$ and $B^2_u \subseteq A_u \mod D_0,$
\sn
\item[$\bullet$]  $B^2_u \in D,$
\sn
\item[$\bullet$]  $\langle B^2_u:u \in [\mu]^{< \theta}\rangle$ is
  multiplicative modulo 
  $D_0$ (see \ref{g13}).
\end{enumerate}

By Definition \ref{g13}(1) applied to $\langle B^2_u:u \in 
[\mu]^{<\theta}\rangle$ recalling clause (c) of the
assumption of the claim, we can find $\bar B = \langle B_u:u \in
[\mu]^{< \theta}\rangle$ which is a multiplicative 
refinement of $\bar B^2$ 
and is multiplicative, 
and $ B_u \in D $ because $ B_n = B^2_n $ modulo 
$D_0  \subseteq D $ and $ B^2_u \in D$.

So we are done for the ``if" direction.

\underline{The ``only if" direction}:

So we are assuming $D$ is a $(\mu,\theta,\varepsilon!,\Delta,T)$-good
filter on $I$ and we have to prove $D_1$ is
$(\mu,\theta,\varepsilon!,\Delta,T)$-moral.

So let $\bar{\mathbf a}$ be a $D_1-(\mu,\theta,\varepsilon!,\Delta,T)$-moral
problem (on $\gB_1$), we have to find a solution. 
For $u \in 
[\mu]^{< \theta}$ choose $A^1_u \subseteq I$ such that $\mathbf
j(A^1_u) = \mathbf a_u$,
so $ A^1_u \in D $ (by clause (e)) and 
$ u \subseteq v \in [\mu ]^{< \theta }\rightarrow 
A^1_u \subseteq A^1_v $  modulo $ D_0 $.
Now 
by \ref{g13}, i.e. clause (b) of the assumption
of the claim there is $\bar A^2 = \langle A^2_u:u \in [\mu]^{<
\theta}\rangle$ such that $A^2_u \subseteq A^1_u,A^2_u = A^1_u \mod D_0$
hence $A^2_u \in D$ and $\bar{ A}^2$ is $\subseteq$-decreasing 
[Why?  Because $\bar {A}^1$ is $\subseteq$-decreasing 
modulo $ D_0 $  
as
$\bar{\mathbf a}$ is 
decreasing hence $\bar A^2$ is
$\subseteq$-decreasing.] 

As $D$ is $(\mu,\theta,\varepsilon !,\Delta,T)$-good filter on $I$
there is a $D$-multiplicative refinement $\langle B^2_u:u \in [u]^{<
\theta}\rangle$
of $ \langle A^2_u : u \in [\mu ]^{< \theta }\rangle $. 
Let $\mathbf b_u = \mathbf j(B^2_u)$, now $\langle \mathbf
b_u:u \in [\mu]^{< \theta}\rangle$ is as required.

2) Similarly.
\end{PROOF}

\begin{claim}
\label{g21} 
Let $ D$ be a filter on $ I $.

1) $D$ is $(\mu,\theta)$-excellent implies $D$ is $(\mu,\theta)$-good, 
see \ref{g13}(3).

2) $D$ is $(\mu,\theta)$-good implies $D$ is
   $(\mu,\theta,\varepsilon,\Delta,T)$-moral. 
\end{claim}

\begin{PROOF}{\ref{g21}}
1) So let $\bar{\mathbf a} = \langle \mathbf a_u:u \in [\mu]^{<
   \theta}\rangle$ be a $D$-problem and we should find a
   $D-(\mu,\theta)$-solution $\bar{\mathbf b}$ below $\bar{\mathbf a}$.
   As $D$ is $(\mu,\theta)$-excellent we apply this to $\bar{\mathbf a}$
   and $\bar{\mathbf b}$ as in \ref{g13}(2).  Easily it is as required.

2) Just read the definitions: there are fewer problems.
\end{PROOF}

\begin{remark}
\label{g40}
We may wonder, e.g. in \ref{g10}(1): can we remove 
the regularity demand on
the filter $D$ from clause (A) to clause (B)?  The answer is yes for most $T$'s.
\end{remark}

\begin{claim}
\label{g43}
The filter $D$ is $(\mu ,\theta)$-regular \when \,:

\begin{enumerate}
\item[$(A)$]  $(a) \quad \gB = \cP(I)$
\sn
\item[${{}}$]  $(b) \quad D$ is a $\theta $-complete ultra-filter on
  $\gB,$
\sn
\item[${{}}$]  $(c) \quad $ the demand 
$\theta > \varepsilon$, is natural but not actually required,
\sn
\item[${{}}$]  $(d) \quad T$ is a complete
  $\bbL_{\theta,\theta}(\tau)$-theory, e.g. $T =
  \Th_{\bbL_{\theta,\theta}}(M),M$ a $\theta$-saturated

\hskip25pt  model (note that $T = T^{[\theta]}_0$ where $T_0 =
  \Th_{\bbL_{\aleph_0,\aleph_0}}(M)$, 

\hskip25pt  i.e. $T$ is determined by $T_0$ and $\theta$),
\sn
\item[$(B)$]  $T$ has a model $M$ and $p = \{\varphi_\alpha(\bar
  x_{[\varepsilon]},\bar b_\alpha):\alpha < \mu\},\varphi_\alpha(\bar
  x_{[\varepsilon]},\bar y_\alpha) \in \bbL_{\theta,\theta},\bar
  b_\alpha \in {}^{\ell g(\bar y_\alpha)}M$ 
    satisfying:  
  
  for every $q
  \subseteq p,$
\sn
\item[${{}}$]  $\quad \bullet \quad q$ is realized in $M$ iff $|q| <
  \theta.$
\sn
\item[$(C)$]  if $M_s$ is a model of $T$ for $s \in I$ then
  $\prod\limits_{s} M_s/D$ is $(\mu^+,\theta,\varepsilon !,\Delta)$-
  saturated.
\end{enumerate}
\end{claim}

\begin{PROOF}{\ref{g43}}
    Should be clear.
\end{PROOF}

\newpage

\section {A counterexample} \label{4}


In \S2 we generalize \cite[Ch.VI,~2.6]{Sh:c} to filters, using the class of relevant $\RSP$'s $\mathbf r$ being  closed under reduced powers (being a Horn class, see \ref{c7}). Can we generalize the result of
Malliaris-Shelah \cite{Sh:998}? 
Here we give a counter-example.

For this, we have to find:

\begin{enumerate}
    \item[$(*)_1$]  $D$ a filter of $\lambda$ such that the partial order $N_1 = (\bbQ,<)^\lambda/D$ satisfies $\gp^*(N_1) =  \kappa_1 + \kappa_2 < \mu^+ \le \gp^*_{\sym}(N_1),\kappa_1 \ne
    \kappa_2,(\kappa_1,\kappa_2) \in \cC(N_1)$, so in fact $N_1$ has no $(\theta_1,\theta_2)$-cut when $\theta_1 = \cf(\theta_1) = \theta_2 \le \mu$ and when $\theta_\ell \ge \mu^+ \wedge \theta_{3-\ell} \in \{0,1\},$
    \sn
    \item[$(*)_2$]  preferably: $\lambda = \mu$
    \sn
    \item[$(*)_3$]  or at least for some dense linear order $M_0$ there is a complete Boolean algebra $\gB$ and a filter $D$ on $\gB$ such that
    $N_0 = M^{\gB}_0/D$ is as above.
\end{enumerate}

We presently deal with the (main) case $\theta=\aleph_0$  and carry this out.  It
seems reasonable that we can prove, e.g. 
$ T_{\ceq} \ntriangleleft_{\rp} T_{\ord}$ but we have not arrived to it; see \cite{Sh:1164} on $ T_{\ceq}$ and  \cite{Sh:457} on the closely related $T_{\feq} $. Later 
we hope 
to say more. Clearly, we can control the set of non-symmetric pre-cuts.

\begin{convention}\label{b0}
    $T_{\ord}$ is the first order theory of $(\bbQ,\le)$, see \ref{b7}(1)(d). 
\end{convention}

\begin{definition}\label{b2}
    Let $\kappa$ be a regular
        uncountable cardinal. 

    1) Let $K^{\ba}_\kappa$ be the class of $\mathbf m$ such that:
    
    \begin{enumerate}
        \item[$(a)$]  $\mathbf m = (\gB,D) = (\gB_{\mathbf m},D_{\mathbf m}),$
        \sn
        \item[$(b)$]  $\gB$ is a complete Boolean algebra satisfying the
          $\kappa$-c.c.,
        \sn
        \item[$(c)$]  $D$ is a filter on $\gB$.
    \end{enumerate}
    
    2) Let $\le^{\ba}_\kappa$ be the following two-place relation on
    $K^{\ba}_\kappa:\mathbf m \le^{\ba}_\kappa \mathbf n$ iff:
    
    \begin{enumerate}
        \item[$(a)$]  $\mathbf m,\mathbf n \in K^{\ba}_\kappa,$
        \sn
        \item[$(b)$]  $\gB_{\mathbf m} \lessdot \gB_{\mathbf n},$
        \sn
        \item[$(c)$]  $D_{\mathbf m} = D_{\mathbf n} \cap \gB_{\mathbf m}$.
    \end{enumerate}
    
    3) Let $S^{\ba}_\kappa$ be the class of $\le^{\ba}_\kappa$-increasing
    continuous sequences $\bar{\mathbf m}$ which means:
    
    \begin{enumerate}
    \item[$(a)$]  $\bar{\mathbf m} = \langle \mathbf m_\alpha:\alpha < \ell g(\bar{\mathbf m})\rangle,$ 
    \sn
    \item[$(b)$]  $\mathbf m_\alpha \in K^{\ba}_\kappa,$
    \sn
    \item[$(c)$]  if $\alpha < \beta < \ell g(\bar{\mathbf{m} })$ then $\mathbf m_\alpha \le^{\ba}_\kappa \mathbf m_\beta,$
    \sn
    \item[$(d)$]  if $\beta < \ell g(\bar{\mathbf m})$ is a limit ordinal
      then:
    \sn
    \begin{enumerate}
        \item[$(\alpha)$]  $\gB_{\mathbf m_\beta }$ is the completion of $\cup\{\gB_{\mathbf m_\alpha}:\alpha < \beta\},$
        \sn
        \item[$(\beta)$]  $D_{\mathbf m_\beta}$ is generated (as a filter) by
          $\cup\{D_{\mathbf m_\alpha}:\alpha < \beta\}$.
    \end{enumerate}
    \end{enumerate}
    
    4) If $\kappa = \aleph_1$ we may write $K^1_{\ba},\le^1_{\ba},S^1_{\ba}$,
    and if $ \kappa = \infty $ we may write 
    $K^2_{\ba},\le^2_{\ba},S^2_{\ba}$
    or $K^{\ba}_\infty ,\le^{\ba}_\infty ,S^{\ba}_\infty $.

    5) We say $\mathbf m$ is of cardinality $ \lambda $ when $\gB_{\mathbf m}$ is  of cardinality 
    $\lambda$.
\end{definition}

\begin{claim}
\label{b4}
1) For every $\lambda$ there is $\mathbf m \in K^{\ba}_\kappa$ of
cardinality $\lambda^{<  \kappa }$.

2) $\le^{\ba}_\kappa$ is a partial order on $K^{\ba}_\kappa$.

3) If $\bar{\mathbf m} = \langle \mathbf m_\alpha:\alpha < \delta\rangle$
is a $\le^{\ba}_\kappa$-increasing continuous sequence, \then \, for
some $\mathbf m_\delta$, the sequence $\bar{\mathbf m} \char 94 \langle
\mathbf m_\delta\rangle$ is $\le^{\ba}_\kappa$-increasing continuous.
\end{claim}

\begin{PROOF}{\ref{b4}}
1) E.g. $\gB_{\mathbf m}$ is the completion of a free Boolean algebra
generated by $\lambda^{< \kappa }$ elements.

2) Easy.

3) If $\cf(\delta) \ge \kappa$, then $\gB_{\mathbf m_\delta} =
\bigcup\limits_{\alpha < \delta} \gB_{\mathbf m_\alpha}$, if
$\cf(\delta) < \kappa$ it is the (pendantically $a$) completion of the
union.  $D_{\mathbf m_\delta}$ is the filter generated by
$\cup\{D_{\mathbf m_\alpha}:\alpha < \delta\}$.  Classically 
$\kappa$-c.c. is preserved.
\end{PROOF}

\begin{definition}
\label{b7}
Let $\mathbf m \in K^2_{\ba}$ and $ \kappa _1, \kappa _2$ are 
(infinite) regular cardinals.

1) We say $\bar{\mathbf a}$ is a $T_{\ord}-(\kappa_1,\kappa_2)$-moral
   problem in $\mathbf m$ \when \,:

\begin{enumerate}
\item[$(a)$]  $\mathbf m \in K^2_{\ba}$, (actually already assumed),
\sn
\item[$(b)$]  $I = I(\kappa_1,\kappa_2)$ is the linear order $I_1 +
  I_2$ where:
\sn
\item[${{}}$]  $\bullet \quad I_1 = I_1(\kappa_1) = 
(\{1\} \times \kappa_1)$,
\sn
\item[${{}}$]  $\bullet \quad I_2 = I_2(\kappa_2) = (\{2\} \times \kappa^*_2),$
\sn
\item[$(c)$]  $\bar{\mathbf a} = \langle \mathbf a_{s,t}:s
  <_{I(\kappa_1,\kappa_2)} t\rangle$ is a sequence of members of $D_{\mathbf m},$
\sn
\item[$(d)$]  if $u \subseteq I$ is finite, $\mathbf t:u \times u
  \rightarrow \{0,1\}$ and $\cap\{\mathbf a^{\iif(\mathbf
  t(s,t))}_{s,t}:s,t \in u\} > 0_{\mathbf{m} },$ \then \, there is a function 
$f:u \rightarrow \{0,\dotsc,|u|-1\}$ such that:
\sn
\item[${{}}$]  $\bullet \quad$ if $s,t \in u$ then $\mathbf t(s,t) = 1$
  \Iff \, $f(s) \le f(t).$
\sn
\item[$(e)$]   hence $s_1 <_I s_2 <_I s_2 \Rightarrow
\mathbf a_{s_1,s_2} \cap \mathbf a_{s_2,s_3} \le \mathbf a_{s_1,s_3}$ and we
stipulate $\mathbf a_{s,s} = 1_{\gB_\mathbf{m} },\mathbf a_{t,s} = \mathbf a_{s,t}$ when
$s <_I t$.
\end{enumerate}

2) We say $\bar{\mathbf b}$ is a solution of $\bar{\mathbf a}$ in $\mathbf
m$ where $\bar{\mathbf a}$ is as above \when \,:

\begin{enumerate}
\item[$(a)$]  $\bar{\mathbf b} = \langle \mathbf b_s:s \in I\rangle,$
\sn
\item[$(b)$]  $\mathbf b_s \in D_{\mathbf m},$
\sn
\item[$(c)$]  if $s_1 \in I_1,s_2 \in I_2$ then $\mathbf b_{s_1} \cap  
\mathbf b_{s_2} \le \mathbf a_{s_1,s_2}$.
\end{enumerate}
\end{definition} 

\begin{definition} \label{b8}
1) For $ \iota =1,2 $ let  
$\mathbf{S}_\iota $ be  the class of tuples 
$\mathbf s = (I,D_0,\mathbf j,\gB,D_1,D)$
such that:

\begin{enumerate}
\item[$(a)$]  $\mathbf j$ is a homomorphism from $\cP(I)$ onto
  the complete Boolean  algebra  $\gB,$
\sn
\item[$(b)$]   $D_1$ is a filter on $\gB,$
\sn
\item[$(c)$]   $D_0 = \{A \subseteq I:\mathbf j(A) =1_{\gB}\}$ (or see
\S3),
\sn
\item[$(d)$]   $D = \{A \subseteq I:\mathbf j(A) \in D_1\}$,
\sn 
\item[(e)] the pair $ (\mathfrak{B} , D ) $ belongs to $ K^\iota _{\ba}$.
\end{enumerate}

2) For $\mathbf s \in \mathbf S$ let $\mathbf m_{\mathbf s} = (\gB_{\mathbf s},D_
\mathbf{s}  
)$.

3) We say $\mathbf s \in \mathbf S$ is
$(\mu,\theta)$-excellent (if $\theta = \aleph_0$ may omit 
it)  
\when \, 
$D_0$ is an excellent filter on $I$, see Definition \ref{g13}(2).

4) We say $\mathbf s \in \mathbf S$ is $(\mu,\theta)$-regular (if $\theta =
   \aleph_0$ we may omit $\theta$) \when \, $D_0$ is a
   $(\mu,\theta)$-regular filter.

5) Let $\mathbf{S}^\iota _{\mu,\theta}$ be the class of $(\mu,\theta)$-excellent
$(\mu,\theta)$-regular $\mathbf s \in \mathbf{S}_\iota $; we may omit $\theta$ if
$\theta = \aleph_0$.

6) Let $ \mathbf{S}_{\mu, \theta, \kappa }$ be the class
of $ \mathbf{s} \in \mathbf{S} ^2_{\mu, \theta }$ such that 
$ \mathfrak{B} _ \mathbf{s} $ satisfies the 
$ \kappa $-c.c.

%
\end{definition}

\begin{claim}
\label{b10}
1) Assume 
$\mathbf m = (\gB,D) \in K_{\ba}$ and $\kappa_1,\kappa_2$ are
infinite and regular cardinals.  
Then for some $M
\in \Mod_{T_{\ord}},M^{\gB}/D$ has a $(\kappa_1,\kappa_2)$-pre-cut \Iff \,
some $T_{\ord}-(\kappa_1,\kappa_2)$-moral problem in $\mathbf m$ has no
solution.

2) Let $ \mu \ge {\aleph_0} = \theta $.
If $\mathbf s \in \mathbf S_{\mu,\theta}$ so is $\mu$-excellent 
and $\mu$-regular and $\kappa_1,\kappa_2 \ge \aleph_0$ are 
regular and $\kappa_1 + \kappa_2 \le \mu$ \then \, the
   following conditions are equivalent:

\begin{enumerate}
\item[$(a)$]  for some linear order $M,M^{I(\mathbf s)}/D_{\mathbf s}$ has
a $(\kappa_1,\kappa_2)$-pre-cut,
\sn
\item[$(b)$]  for every infinite linear order, $M^{I(\mathbf s)}/D_{\mathbf
s}$ has a $(\kappa_1,\kappa_2)$-pre-cut,
\sn
\item[$(c)$]  not  every $T_{\ord}-(\kappa_1,\kappa_2)$-moral problem in
$\mathbf m_{\mathbf s}$ has a solution.
\end{enumerate}
\end{claim}

\begin{PROOF}{\ref{b10}}
As in    
    the proof of 
\ref{g10}(1), 
relying on Def \ref{b7} instead of  Def \ref{g3}; 
recalling

\begin{enumerate}
\item[$\boxplus$]  if $M^\iota_s$ for $s \in I,\iota \in \{1,2\}$ are
  $\tau$-models, $|\tau| \le \mu,D$ a $\mu$-regular filter on $I$ and
  $M^1_s,M^2_s$ are elementarily equivalent, \then \, $N_1 =
  \prod\limits_{s \in I} M^1_s/D,N_2 = \prod\limits_{s \in I} M^2_s/D$
  are $\bbL_{\mu^+,\mu^+}$-equivalent (and more, see Kennedy-Shelah
  \cite{Sh:769}, \cite{Sh:852} and 
  Kennedy-Shelah-Vaananen \cite{Sh:912} on the subject).
\end{enumerate}
\end{PROOF}

\begin{observation}
\label{b12}
Assume 
$ \mathbf{m} \in K^2_{\ba}$ and 
$\bar{\mathbf a}$ is a $T_{\ord}-(\kappa_1,\kappa_1)$-moral
problem for $\mathbf m$ so 
(see \ref{b8}(5)) $I_\ell = I_\ell(\kappa_\ell)$
for $\ell=1,2$.


1) If $I'_1 \subseteq I_1$ is cofinal in $I_1$ and $I'_2 \subseteq I_2$
is co-initial in $I_2$ \then \, $\bar{\mathbf a}$ has a solution in $\mathbf m$
\Iff \, $\bar{\mathbf a}' = \bar{\mathbf a} \rest (I'_1 + I'_2) = \langle \mathbf
   a_{s,t}:s <_I t$ and $s,t \in I'_1 + I'_2\rangle$ has 
a solution in $\mathbf m$.

1A) Also, above, if $\bar{\mathbf b}$ is a solution of $\bar{\mathbf a}$ in
$\mathbf m$, \then \, $\bar{\mathbf b} \rest (I'_1 + I'_2)$ is a solution of
$\bar{\mathbf a}'$.

1B) Also above, if $\bar{\mathbf b}'$ is a solution of $\bar{\mathbf a}'$,
\then \, $\bar{\mathbf b}$ is a solution of $\bar{\mathbf a}$ \when \,:

\begin{enumerate}
\item[$(a)$]  if $s \in I_1$ and $t \in I'_1$ is minimal such that $s
  \le_I t$ then $\mathbf b_s = \mathbf b'_t \cap \mathbf a_{s,t}$ if $s <_I
  t$ and $\mathbf b_s = \mathbf b'_t$ if $s=t,$
\sn
\item[$(b)$]  like (a) replacing $I_1,I'_1,s <_I t,\mathbf a_{s,t}$ by
  $I_2,I'_2,t \le_I s,\mathbf a_{t,s}$.
\end{enumerate}

2) If $\bar{\mathbf b}$ is a solution of $\bar{\mathbf a}$ in $\mathbf m$
and $\mathbf b'_s \in D \wedge \mathbf b'_s \le \mathbf b_s$ for $s \in I_1
+ I_2$ then $\langle \mathbf b'_s:s \in I\rangle$ is a 
solution of $\bar{\mathbf a}$ for $\mathbf m$.
\end{observation}

\begin{PROOF}{\ref{b12}}
1) Easy using the proofs of \ref{g10}, \ref{b10} or using (1A), (1B).

1A), 1B), 2)  Check.
\end{PROOF}

A key point in the inductive construction is:

\begin{claim}
\label{b13}
There is no solution to $\bar{\mathbf a}$ in $\mathbf m_\delta$ \when \,:

\begin{enumerate}
\item[$(a)$]  $\bar{\mathbf m} = \langle \mathbf m_\alpha:\alpha \le
  \delta\rangle \in S^2_{\ba},$
\sn
\item[$(b)$]  $\bar{\mathbf a}$ is a
  $T_{\ord}-(\kappa_1,\kappa_2)$-moral problem in $\mathbf m_0,$
\sn
\item[$(c)$]  if $\alpha < \delta$ then $\bar{\mathbf a}$ has no
  solution in $\mathbf m_\alpha,$
\sn
\item[$(d)$]  $\cf(\delta) \ne \kappa_1$ or $\cf(\delta) \ne
  \kappa_2$.
\end{enumerate}
\end{claim}

\begin{PROOF}{\ref{b13}}
Let $\mathbf m_\gamma = (\gB_\gamma,D_\gamma)$ for $\gamma \le \delta$;
by symmetry \wilog \, $\cf(\delta) 
\ne \kappa_1$ and toward contradiction assume $\bar{\mathbf b} = 
\langle \mathbf b_s:s \in I_1 + I_2\rangle$ is a solution of 
$\bar{\mathbf a}$ in $\mathbf m_\delta$.

Hence $\mathbf b_s \in D$.
Now $D_\delta$ is not necessarily equal to
$\bigcup\limits_{\gamma < \delta} D_\delta$ but recalling
\ref{b2}(3)(d)$(\beta)$ and $\langle D_\gamma:\gamma < \delta\rangle$
being increasing, clearly, every member of $D_\delta$ is above some
member of $\bigcup\limits_{\gamma < \delta} D_\gamma$.

So by Observation \ref{b12}(2), \wilog \, $s \in I_1 + I_2
\Rightarrow \mathbf b_s \in \bigcup\limits_{\gamma < \delta} D_\gamma
\subseteq \bigcup\limits_{\gamma < \delta} \gB_\gamma$.

As $\cf(\delta) \ne \kappa_1$, for some $\gamma < \delta$ we have $\kappa_1 =
\sup\{\alpha < \kappa_1:\mathbf b_{(1,\alpha)} \in \gB_\gamma\}$,
i.e. $\{s \in I_1:\mathbf b_s \in \gB_\gamma\}$ is co-final in $I_1$.
So by \ref{b12}(1) \wilog \, 

\begin{enumerate}
\item[$(a)$]  $s \in I_1 \Rightarrow \mathbf b_s \in \gB_\gamma$.
\end{enumerate}

As $D_\gamma = D_\delta \cap \gB_\gamma$ by \ref{b2}(2)(c) clearly

\begin{enumerate}
\item[$(b)$]  $s \in I_1 \Rightarrow \mathbf b_s \in D_\gamma$.
\end{enumerate}

For $t \in I_2$ let $\mathbf b'_t = \min\{\mathbf b \in 
\gB_\gamma:\gB_\delta \models \mathbf b_t \le \mathbf b\}$, well defined
  because $\gB_\gamma$ is complete.

Now

\begin{enumerate}
\item[$(c)$]  $\mathbf b'_t \in D_\gamma$ for $t \in I_2$.
\end{enumerate}

[Why?  Clearly $\mathbf b_t \in \gB_\delta$ as 
$\bar{\mathbf b}$ is a solution of $\bar{\mathbf a}$ in $\mathbf m_\delta$
and $\mathbf b_t \le \mathbf b'_t,\mathbf b'_t \in \gB_\gamma$ by its choice.
Also $\mathbf b'_t \in D_\delta$ because $\mathbf b_t \le \mathbf b'_t \wedge
\mathbf b_t \in D_\delta$ and $D_\delta$ is a filter on $\gB_\delta$ and
lastly $\mathbf b'_t \in D_\gamma$ as $D_\gamma = D_\delta \cap \gB_\gamma$.]

\begin{enumerate}  
\item[$(d)$]  if $s \in I_1,t \in I_2$ then $\mathbf b_s \cap \mathbf b'_t
  \le \mathbf a_{s,t}$.
\end{enumerate}

[Why?  Note $\gB_\delta \models ``\mathbf b_s \cap \mathbf b_t
  \le \mathbf a_{s,t}"$ because $\bar{\mathbf b}$ is a solution of $\mathbf
a$ in $\gB_\delta$ hence
  $\mathbf b_t \le \mathbf a_{s,t} \cup (1-\mathbf b_s)$ and the later $\in
  \gB_\gamma$.  So by the choice of $\mathbf b'_t,\mathbf b'_t 
\le \mathbf a_{s,t} \cup (1- \mathbf b_s)$ hence $\mathbf b_s \cap
  \mathbf b'_t \le \mathbf a_{s,t}$.]

\begin{enumerate}
\item[$(e)$]  $\langle \mathbf b_s:s \in I_1 \rangle \char 94 \langle \mathbf
  b'_t:t \in I_2 \rangle$ solves $\bar{\mathbf a}$ in $\gB_\gamma$.
\end{enumerate}

[Why?  By (a) + (b) + (c) + (d).]

But this contradicts an assumption.
\end{PROOF}

\begin{definition}
\label{b17}
Assume $\mathbf m \in K^2_{\ba}$ and $\bar{\mathbf a}$ is a
$
T_{\ord}-  
(\kappa_1,\kappa_2)$-moral problem in $\mathbf m$.  We say $\mathbf n$ is a
simple $\bar{\mathbf a}$-solving extension of $\mathbf m$ \when \,:

\begin{enumerate}
\item[$(a)$]  $\gB_{\mathbf n}$ is the completion of $\gB^o_{\mathbf n}$
  where:
\sn
\item[$(b)$]  $\gB^o_{\mathbf n}$ is
the Boolean algebra 
generated by $\gB_{\mathbf m}
  \cup \{y_s:s \in I(\kappa_1,\kappa_2)\}$ freely except the 
equations which holds in $\gB_{\mathbf m}$ and 
$\Gamma_{\bar{\mathbf a}} = \{y_{s_1} \cap y_{s_2}
  \le \mathbf a_{s_1,s_2}:s_1 \in I_1(\kappa_1)$ and $s_2 \in I_2(\kappa_2)\},$
\sn
\item[$(c)$]  $D_{\mathbf n}$ is the filter on $\gB_{\mathbf n}$ 
generated by $D_{\mathbf m} \cup \{y_s:s \in I(\kappa_1,\kappa_2)\}$.
\end{enumerate}
\end{definition}


\begin{claim}
\label{b20}
Assume $\bar{\mathbf a}$ is a $T_{\ord}-(\kappa_1,\kappa_2)$-moral
problem in $\mathbf m \in K^{\ba}_\kappa$ and\footnote{It seems that
  $\min\{\kappa_1,\kappa_2\} < \kappa$ suffice; the only difference in
  the proof is in proving $(*)_5$.}
$\kappa = \cf(\kappa) > \kappa_1 + \kappa_2$.

1) There is $\mathbf n \in K^{\ba}_ \kappa $ which is a simple 
$\bar{\mathbf a}$-solving extension of $\mathbf m$, unique up to
isomorphism over $\gB_{\mathbf m}$.

2) Above $\mathbf m \le^{\ba}_\kappa \mathbf n$ (so $\mathbf n \in
   K^{\ba}_\kappa$).

3) If $\bar{\mathbf a}^*$ is a $T_{\ord}-(\theta_1,\theta_2)$-moral
   problem of $\mathbf m$ with no solution in $\mathbf m$ and $\theta_1
   \notin \{\kappa_1,\kappa_2\}$ or $\theta_2 \notin
   \{\kappa_1,\kappa_2\}$ \then \, $\bar{\mathbf a}^*$ has 
no solution in $\mathbf n$.
\end{claim}

\begin{PROOF}{\ref{b20}}
1) As above let $I_\ell = I_\ell(\kappa_\ell)$ for $\ell=1,2$ and $I=I_1 + I_2$.

First,

\begin{enumerate}
\item[$(*)_1$]  the set of equations $\Gamma_{\bar{\mathbf a}}$ is
  finitely satisfiable in $\gB_{\mathbf m}$.
\end{enumerate}

Why?  
We  shall prove two 
stronger statements (each implying $(*)_1$).

\begin{enumerate}
\item[$(*)_{1.1}$]  if $t_1 \in I_1$ then we can find $\langle \mathbf b'_s:s
  \in I\rangle \in {}^I \gB$ such that:
\sn
\begin{enumerate}
\item[$(a)$]  $\mathbf b'_s \in D_{\mathbf m} \subseteq \gB_m$
if $(s \le _{I_1}
  t_1) \vee (s \in I_2),$
\sn
\item[$(b)$]  if $s_1 \in I_1,s_2 \in I_2$ then $\mathbf b'_{s_1} \cap
  \mathbf b'_{s_2} \le \mathbf a_{s_1,s_2}$.
\end{enumerate}
\end{enumerate}

[Why?  Let $\mathbf b'_s$ be:

\begin{enumerate}
\item[$\bullet$]  $\mathbf a_{s,t_1}$ if $s \le _I t_1$ (so $s \in I_1$),
\sn
\item[$\bullet$]  $\mathbf a_{t_1,s}$ if $s \in I_2,$
\sn
\item[$\bullet$]  $0_{\gB}$ if $t_1 <_I s \in I_1$.
\end{enumerate}

Now clause (a) is obvious 
(recalling $ \mathbf{a} _{t_1, t_1}= 1_{\mathfrak{B} _\mathbf{m} }$
and as for 
clause (b), let $s_1 \in I_1,s_2 \in
I_2$, now if $t_1 \le_I s_1 \in I_1$ then $\mathbf b'_{s_1} \cap \mathbf
b'_{s_2} = 0_{\mathfrak{B} _\mathbf{m} } \cap \mathbf b'_{s_2} 
= 0_{\mathfrak{B} _\mathbf{m} } \le \mathbf a_{s_1,s_2}$ and if $s_1
<_I t_1$ then $\mathbf b'_{s_1} \cap \mathbf b'_{s_2} = \mathbf a_{s_1,t_1}
\cap \mathbf a_{t_1,s_2}$ which is $\le \mathbf a_{s_1,s_2}$ by
\ref{b7}(1)(d),(e).]

\begin{enumerate}
\item[$(*)_{1.2}$]  if $t_2 \in I_2$ then we can find $\langle \mathbf b'_s:s
  \in I\rangle \in {}^I \gB$ such that:
\sn
\begin{enumerate}
\item[$(a)$]  $\mathbf b'_s \in D_{\mathbf m} \subseteq \gB_{\mathbf m}$ if
  $s \in I_1$ or $t_2 \le _{I_2} s,$
\sn
\item[$(b)$]  if $s_1 \in I_2,s_2 \in I_2$ \then \, $\mathbf b'_{s_1}
  \cap \mathbf b'_{s_2} \le \mathbf a_{s_1,s_2}$.
\end{enumerate}
\end{enumerate}

[Why?  Similarly.]

Now $(*)_1$ is easy: if $\Gamma' \subseteq \Gamma_{\bar{\mathbf a}}$ is
finite let $t_* \in I_1$ be such that: if $t \in I_1$ and $y_t$
appears in $\Gamma'$ then $t \le _I t_*$. 
Choose $\langle b'_s:s \in
I\rangle$ as in $(*)_{1.1}$ for $t_*$ and let $h$ be the function $y_s
\mapsto \mathbf b'_s$ for $s \in I$.  Now think, 
so $ ( * ) _1$ holds indeed.

Clearly it follows by $(*)_1$ that

\begin{enumerate}
\item[$(*)_2$]  $(a) \quad$ there is a Boolean algebra $\gB^o_{\mathbf
  n}$ extending $\gB_{\mathbf m}$ as described in 

\hskip25pt  clause (b) of Definition \ref{b17},
\sn
\item[${{}}$]  $(b) \quad$ there is a Boolean algebra $\gB_{\mathbf n}$
as described in (a) of Definition \ref{b17}: 

\hskip25pt  the completion of $\gB^0_{\mathbf n},$
\sn
\item[${{}}$]  $(c) \quad D_{\mathbf n}$ is 
chosen as the filter on $\gB_{\mathbf n}$
  generated by $D_{\mathbf m} \cup \{y_s:s \in I\}$ 

\hskip25pt  it  
    satisfies   
$D_{\mathbf m} = D_{\mathbf n} \cap \gB_{\mathbf m}$,
in particular $0_{\mathfrak{B} _\mathbf{m} } \notin D_{\mathbf n}$,
\sn
\item[${{}}$]  $(d) \quad \gB_{\mathbf n}$ satisfies the $\kappa$-c.c.,
\sn
\item[${{}}$]  $ (e )\quad D_ \mathbf{n}  $  is generated (as a
 filter) by $D_\mathbf{n}  \cap \mathfrak{B}^o  _ \mathbf{n}.$
\end{enumerate}

[Why?  Clauses (a), (b) follows by $(*)_1$ and for clauses (c),(d) 
see $(*)_4$ and $(*)_5$ in the proof of (2), respectively; in
particular $0_{\mathfrak{B} \mathbf{m} } \notin D_{\mathbf n}$
    and clause (e) holds by clause (c).]  

Together we have $\mathbf n = (\gB_{\mathbf n},D_{\mathbf n}) \in K^2_{\ba}$,
as for $\mathbf m \le_{\ba} \mathbf n$, see part (2).

2) Now (by part (1) we have $\gB_{\mathbf m} \subseteq \gB_{\mathbf n}$, but 
we shall show that 
moreover) 

\begin{enumerate}
\item[$(*)_3$]  $\gB_{\mathbf m} \lessdot \gB_{\mathbf n}$.
\end{enumerate}

[Why?  If not, then some $\mathbf d \in \gB^+_{\mathbf n}$ 
 is disjoint to $\mathbf b$ for a dense subset of $\mathbf b \in \
 \gB_{\mathbf m}^+ $. 
 Let $\mathbf d = \sigma(y_{s_0},\dotsc,y_{s_{n-1}},\bar c)$
where $\sigma$ is a Boolean term, $s_0 <_I \ldots <_I s_{n-1}$
and $\bar c$ is from $\gB_{\mathbf m}$.  We may replace $\mathbf d$ by 
any $\mathbf d' \in \gB^+_\mathbf{n} $ 
satisfying $\mathbf d' \le_{\gB} \mathbf d$.
Hence \wilog \, $ \mathbf{d} = 
\cap\{y^{\iif(\eta(\ell))}_{s_\ell}:
\ell < n\} \cap c > 0_{\mathbf{n} }$ where 
$c \in \gB_{\mathbf m},\eta(\ell) \in \{0,1\}$ 
for $\ell < n$; also 
\wilog \,
for every $\ell,k < n$ we have $s_\ell \in I_1 \wedge s_k \in I_2
\Rightarrow (c \le \mathbf a_{s_\ell,s_k}) \vee (c
\cap \mathbf a_{s_\ell,s_k} = 0_{\mathfrak{B} _\mathbf{n} })$.

We now define a function $h$ from $\{y_s:s \in I\}$ into $\gB_{\mathbf m}$ 
as follows: $h(y_s) \text{ is}$:

\begin{enumerate}
\item[$\bullet_1$]  $c, \text{ if } \quad s=s_\ell \wedge \eta(\ell)=1,$
\sn
\item[$\bullet_2$]  $ 0_{\mathfrak{B} _\mathbf{m} },$  \text{ if otherwise}.
\end{enumerate}

Now

\begin{enumerate}
\item[$\bullet_3$]  if $t_1 \in I_1,t_2 \in I_2$ then $\gB_{\mathbf m}
  \models ``h(y_{t_1}) \cap h(y_{t_2}) \le \mathbf a_{t_1,t_2}"$.
\end{enumerate}

[Why?  If $h(t_1)=0_{\mathfrak{B} _ \mathbf{m} } \vee h(t_2)=
0_{\mathfrak{B} _\mathbf{m} }$ this is obvious, otherwise for some
$\ell(1) < \ell(2) < n$ we have $t_1 =
s_{\ell(1)},t_2 = s_{\ell(2)}$ and $\eta(\ell(1)) = 1 =
\eta(\ell(2))$.  So 
it suffice 
to prove $c = c \cap c \le \mathbf a_{t_1,t_2}$ but otherwise by
the choice of $c,c \cap \mathbf a_{t_1,t_2}=0$, hence recalling
\ref{b17}(b) we have $\gB_{\mathbf n}
\models ``y_{s_1} \cap y_{s_2} \cap c = 0"$ contradiction to our
current assumption $\gB_{\mathbf n} \models ``d>0$"; so $\bullet_3$
holds indeed.]

By the choice of
$ \Gamma _ { \bar{ \mathbf{a} }} $ and of 
$\gB_{\mathbf n}$ recalling $\gB_{\mathbf m}$ is complete, by 
the choice of $h$ and $\bullet_3$ there is a projection 
$\hat h$ from $\gB_{\mathbf n}$ onto
$\gB_{\mathbf m}$ extending $h$, so clearly $\hat h(d)=c$ and this
implies $c_1 \in \gB_{\mathbf m} \wedge 0 < c_1 \le c \Rightarrow
\gB_{\mathbf n} \models ``c_1 \cap \mathbf{d}  \ge 0_{\mathfrak{B} _
\mathbf{n} }"$ contradicting the choice of
$\mathbf{d} $.  So indeed $(*)_3$ holds.]

\begin{enumerate}
\item[$(*)_4$]  $D_{\mathbf m} = D_{\mathbf n} \cap \gB_{\mathbf m}$.
\end{enumerate}

[Why?  Otherwise there are $c_1 \in D_{\mathbf m},c_2 \in \gB_{\mathbf m}
  \backslash D_{\mathbf m}$ and $s_0 <_I \ldots <_I s_{n-1}$ such that
  $\gB_{\mathbf n} \models ``\bigcap\limits_{\ell < n} y_{s_\ell}  \cap
  c_1 \le c_2"$.  As $\mathbf a_{t_1,t_2} \in D_{\mathbf m}$ for 
$t_1 <_I t_2$, \wilog \, $c_1 \le 
\mathbf a_{s_\ell,s_k}$ for $\ell < k < n, 
s_ {\ell} \in I_1, s _ k \in I_2$.   

Now letting $c = c_1-c_2$ we 
continue as in the proof of $(*)_3$ defining $h,\hat h$ 
and apply the projection $\hat h$ 
to ``$\bigcap\limits_{\ell < n} y_{s_\ell} \cap c_1 \le c_2$".]

\begin{enumerate}
\item[$(*)_5$]  $\gB_{\mathbf n}$ satisfies the $\kappa$-c.c..
\end{enumerate}

[Why?  If not, then there are pairwise disjoint, positive $d_i \in
\gB_{\mathbf n}$ for $i < \kappa$.  So as in the proof of $(*)_3$,
  \wilog \, $d_i = \cap\{y^{\iif(\eta(i,\ell))}_{s(i,\ell)}:\ell <
 n(i)\} \cap c_i$ where $c_i \in \gB_{\mathbf m},\eta(i,\ell) \in
  \{0,1\}$ and $s(i,0) <_I s(i,1) <_I \ldots <_I s(i,n(i)-1)$.
Let $ m(i) \le n(i) $ be such that for every $ {\ell} < n(i)$
we have $ s_ {\ell} \in I_1 $ iff $ {\ell} < m(i)$. 

Again as there, \wilog \, for every $\ell < m(i) \le k<n(i)$ 
we have $(\mathbf
a_{s(i,\ell),s(i,k)} \le c_i) \vee 
(\mathbf a_{s(i,\ell),s(i,k)} \cap c_i=0)$
so $\eta(i,\ell)=1=\eta(i,k) \wedge \ell < m(i) \le  k < n(i)
\Rightarrow c_i \le
\mathbf a_{s(i,\ell),s(i,k)}$.

As $\kappa = \cf(\kappa) > \kappa_1 + \kappa_2$ by an assumption of \ref{b20} 
\wilog \, $n(i) = n, m(i) = m, 
\eta(i,\ell) = 
\eta(\ell)$ and $s(i,\ell) = s_\ell$ for $i <
\kappa,\ell < n$ and as $\gB_{\mathbf m}$ satisfies the 
$\kappa$-c.c. we
can find $i < j < \kappa$ such that $\gB_{\mathbf m} \models ``0 < c_i \cap
c_j"$ and let $c =c_i \cap c_j$ so we   can 
continue as before.]

So together by $(*)_3, (*)_4, (*)_5$ we have $\mathbf m
\le^{\ba}_\kappa \mathbf n \in K^{\ba}_\kappa$ as promised.

3) Let $I^* = I(\theta_1,\theta_2),I^*_1 = I_1(\theta_1),
I^*_2 = I_2(\theta_2)$ and recall $\bar{\mathbf a}^* = \langle 
\mathbf a^*_{s,t}:s <_{I^*} t\rangle$ is a
$T_{\ord}-(\theta_1,\theta_2)$-moral problem in $\mathbf m$.  Toward
a 
contradiction assume that the sequence $\bar{\mathbf b} = \langle \mathbf b_t:t \in
 I^*\rangle$ solve the problem 
 $\bar{\mathbf a}^*$ in $\mathbf n$ so 
 $ \mathbf{b} _ t \in D_ \mathbf{n} $ and 
 let
$\mathbf b_t = \sigma_t(y_{s(t,0)}\dotsc,y_{s(t,n(t)-1)},
c_{t,0},\dotsc,c_{t,m(t)-1})$ with $c_{t,k} \in \gB_{\mathbf m},
s(t,\ell) \in I$ and \wilog \, $s(t,\ell) <_I s(t,\ell +1)$ 
for $\ell < n(t)-1$ so $s(t,k) \in I$ for $k < n(t)$.

The reader may wonder: we have to prove that there is 
no solution in $ \mathfrak{B} _ \mathbf{n} $ , not just
 in $  \mathfrak{B} ^o _ { \mathbf{n} }$ , so how can we use finitary
 terms? The point is that though 
 $  \mathfrak{B} _ \mathbf{n} $ is the 
 completion of $ \mathfrak{B} ^o_\mathbf{n} $,
 the filter $D_ \mathbf{n} $ is generated (as a filter) 
 by $ \mathfrak{B} ^o_\mathbf{n} \cap D_\mathbf{n} $.

By symmetry \wilog,

\begin{enumerate}
\item[$(*)_6$]   $\theta_1 \notin \{\kappa_1,\kappa_2\}$.
\end{enumerate}

Recalling \ref{b12}, we can replace $\mathbf b_t$ by any $\mathbf b'_t \le
\mathbf b_t$ which is from $D_{\mathbf n}$, so as $\bigwedge\limits_{\ell}
y_{s(t,\ell)} \in D_{\mathbf n}$, \wilog \, 
$\ell < n(t) \Rightarrow
\mathbf b_t \le y_{s(t,\ell)}$, so \wilog \, 
\begin{enumerate} 
\item[$(*)_7$] $\mathbf b_t = 
\cap\{y_{s(t,\ell)}:\ell < n(t)\} \cap c_t$ 
for some $c_t \in D_{\mathbf
m}$ recalling $D_{\mathbf m} = D_{\mathbf m} \cap \gB_{\mathbf m}$.
\end{enumerate}

By the $\Delta$-system lemma (recalling \ref{b12}(1)) \wilog,

\begin{enumerate}
\item[$\oplus$]  if $\theta_1 > \aleph_0$ then,
\sn
\begin{enumerate}
\item[$(a)$]  $t \in I^*_1 \Rightarrow n(t) = n(*),$
\sn
\item[$(b)$]  if $t \in I^*_1$ then $s(t,\ell) \in I^*_1
  \Leftrightarrow \ell < \ell(*),$
\sn
\item[$(c)$]  $\big< \langle s(t,\ell):\ell <
  n(*)\rangle:t \in I^*_1 \big>$ is an indiscernible sequence
  in 

\hskip25pt the linear order $I=I(\kappa_1, \kappa_2 )$,
for quantifier free formulas.
\end{enumerate}
\end{enumerate}

But we shall not use $\oplus$.
As $\theta_1 \ne \kappa_1,\kappa_2$, 
by \ref{b12}(1),(1A)
it follows that 
\wilog 
for some
$s^\circ_1,s^\circ_2$ we have:

\begin{enumerate}
\item[$(*)_8$]  $s^\circ_1 \in I_1,s^\circ_2 \in I_2$ and $s(t,\ell)
  \notin [s^\circ_1,s^\circ_2]_I$ for every $t \in I^*_1,\ell < n(t)$.
\end{enumerate}

Again by \ref{b12}(2) \wilog,

\begin{enumerate}
\item[$(*)_9$]  if $t \in I^*_2$ then $\mathbf b_t \le y_{s^\circ_1} \cap 
y_{s^\circ_2}$.
\end{enumerate}

We now define a function $h$ from 
$\{y_s:s \in I\}$ into $\gB_{\mathbf
  n}$, (yes! not $\gB_{\mathbf m}$) by:

\begin{enumerate}
\item[$(*)_{10}$]  $h(y_s)$ is:
\sn
\begin{enumerate}
\item[$\bullet$]  $\mathbf a_{s,s^\circ_1} \cap \mathbf
a_{s^\circ_1,s^\circ_2}$ \If \, $s <_I s^\circ_1,$
\sn
\item[$\bullet$]  $\mathbf a_{s^\circ_1,s} \cap y_s \cap \mathbf
  a_{s,s^\circ_2}$ \If \, $s \in I,s^\circ_1 \le_I s \le_I s^\circ_2,$
\sn
\item[$\bullet$]  $\mathbf a_{s^\circ_1,s^\circ_2} \cap \mathbf a_{s^\circ_2,s}$ 
\If \, $s^\circ_2 <_I s$.
\end{enumerate}
\end{enumerate}

Note

\begin{enumerate}
\item[$(*)_{11}$]  $h(y_s) \in D_{\mathbf n}$ for $s \in I$.
\end{enumerate}

[Why?  Because $\mathbf a_{s,t} \in D_\mathbf{n} $ 
for     $ S \in I_1, t \in I)2$ 
and $y_s \in D_\mathbf{n} $ for
$s \in I$.]

\begin{enumerate}
\item[$(*)_{12}$]  $h(y_{s_1}) \cap h(y_{s_2}) \le \mathbf a_{s_1,s_2}$ for
  $s_1 \in I_1,s_2 \in I_2$.
\end{enumerate}

[Why?  If $s_1,s_2 \in [s^\circ_1,s^\circ_2]_I$ this holds by the
  definition of $\gB_{\mathbf n}$, i.e. as $h(y_{s_1}) \le y_{s_1},
  h(y_{s_2})
  \le y_{s_2}$ and $\gB_{\mathbf n} \models ``y_{s_1} \cap y_{s_2} \le
  \mathbf a_{s_1,s_2}"$.

If $s_1 <_{I^*} s^\circ_1 \wedge s^\circ_2 <_{I^*} s_2$ then
$(*)_{11}$ says:
$\mathbf a_{s_1,s^\circ_1} \cap \mathbf a_{s^\circ_1,s^\circ_2} \cap
\mathbf a_{s^\circ_2,s_2} \le \mathbf a_{s_1,s_2}$ which obviously holds
(as $\bar{\mathbf a}$ is a $T_{\ord}-(\kappa_1, \kappa_2)$-problem 
in $ \mathbf{m} $). 

If $s_1 <_{I^*} s^\circ_1 \wedge s_2 \in [s^\circ_1,s^\circ_2]_{I^*}$
then this means:
$(\mathbf a_{s_1,s^\circ_1} \cap \mathbf a_{s^\circ_1,s^\circ_2}) \cap
(\mathbf a_{s^\circ_1,s_2} \cap y_{s_2} \cap \mathbf a_{s,s^\circ_2}) 
\le \mathbf a_{s_1,s_2}$; but as we have
$\mathbf a_{s_1,s^\circ_1} \cap \mathbf a_{s^\circ_1,s_2} \le 
\mathbf a_{s_1,s_2}$ this holds.

If $s_1 \in [s^\circ_1,s^\circ_2]_{I^*}$ and $s^\circ_2 <_{I^*} s_2$
this means $(\mathbf a_{s^\circ_1,s_1} \cap y_{s_1} \cap \mathbf a_{s_1,s_2}) \cap
(\mathbf a_{s^\circ_1,s^\circ_2} \cap \mathbf a_{s^\circ_2,s}) \le \mathbf
a_{s_1,s_2}$ which holds for similar reasons.  So $(*)_{12}$ holds
indeed.]

By the choice of $\gB^\circ_{\mathbf n}$ and $\gB_{\mathbf n}$ there is a
homomorphism $\hat h$ from $\gB_{\mathbf n}$ into $\gB_{\mathbf n}$, extending
$\id_{\gB_{\mathbf m}}$ and extending $h$. 
Now easily $\hat h(\mathbf b_t) \in D$
for $t \in I^*$ because $\mathbf b_t = 
\cap\{y_{s(t,\ell)}:\ell < n(t)\}
\cap c_t,c_t \in D_{\mathbf m}$ hence 
$\hat h(c_t) = c_t \in D_{\mathbf m}$ and by
$(*)_{10}$ we have $\hat h(y_{s(\ell,t)}) \in D_{\mathbf m}$.

Now $\langle \hat h(\mathbf b_t):t \in I^*\rangle$ still form a
solution of $\bar{\mathbf a}^*$ and by $(*)_7+(* )_8 +
(*)_{10}$ we have 
$t \in I^*_1 \Rightarrow h(\mathbf b_t) \in
\gB_{\mathbf m}$ hence \wilog \,: 

\begin{enumerate}
\item[$(*)_{13}$]  $t \in I^*_1 \Rightarrow \mathbf b_t \in \gB_{\mathbf m}$.
\end{enumerate}

Now define $\mathbf b'_t$ for $t \in I^*$ by: $\mathbf b'_t$ is:

\begin{enumerate}
\item[$\bullet$]  $\mathbf b_t$ if $t \in I^*_1,$
\sn
\item[$\bullet$]  $c_t$ if $t \in I^*_2$.
\end{enumerate}

It suffices to prove that $\langle \mathbf b'_t:t \in I^*\rangle$ solves
$\bar{\mathbf a}^*$ in $\mathbf m$.  Clearly $t \in I^* \Rightarrow \mathbf
b'_t \in D_{\mathbf m}$, so let $t_1 \in I^*_1,t_2 \in I^*_2$.  
We have to prove that $\mathbf b'_{t_1} \cap \mathbf b'_{t_2} \le 
\mathbf a_{t_1,t_2}$ but we know only that
$\mathbf b_{t_1} \cap \mathbf b_{t_2} \le \mathbf a_{t_1,t_2}$ which means
$\mathbf a_{t_1,t_2} \ge \mathbf b'_{t_1}
\cap (\bigcap\limits_{\ell < n(t_2)} y_{s(t_2,\ell)} \cap c_{t_2}) = 
(\mathbf b'_{t_1} \cap \mathbf b'_{t_2}) \cap 
\bigcap\{y_{s(t_2,\ell)}:\ell < n(t_2)\}$.

Let $h_{t_2}$ be a projection from $\gB_{\mathbf n}$ onto $\gB_{\mathbf m}$
such that $h_{t_2}(y_{s(t_2,\ell)}) = c_t$ if $\ell < n(t)$ and
$h_{t_2}(y_s)=0_{\mathfrak{B} \mathbf{m} }$ if $s \in I \backslash
\{s(t_2,\ell):\ell < n(t_2)\}$, as earlier it exists and applying it we
get the desired inequality.
\end{PROOF}

\begin{theorem}
\label{b23}
For any $\lambda$ and regular $\theta_1,\theta_2 \le \lambda$ such
that $\theta_1 + \theta_2 > \aleph_0$ 
there is a regular filter $D$ on $\lambda$ such that:

\begin{enumerate}
\item[$(a)$]   for every dense
linear order $M$, in $M^\lambda/D$ there is a
$(\theta_1,\theta_2)$-pre-cut but no $(\kappa_1,\kappa_2)$-pre-cut when
$\kappa_1,\kappa_2$ are regular $\le \lambda$ and
$\{\theta_1,\theta_2\} \nsubseteq \{\kappa_1,\kappa_2\}$
\sn
\item[$(b)$]  if $M$ is $({}^{\omega >}2,\triangleleft)^\lambda/D$
  then $\gt(M) \ge \lambda^+$.
\end{enumerate}
\end{theorem}

\begin{remark}
\label{b24}
1) Why do we need $\theta_1 + \theta_2 > \aleph_0$?  To prove $(*)_1$.

2) In fact, this demand is necessary, see \ref{b31} below.
\end{remark}

\begin{PROOF}{\ref{b23}}
We  prove clause (a),
which is the main result,  clause (b) holds 
by \ref{b34}.
Let $\kappa = \lambda^+$.

\begin{enumerate}
\item[$(*)_1$]  there are $\mathbf m_0,\mathbf a$ such that:
\sn
\begin{enumerate}
\item[(a)]   $\mathbf m_0 \in K^{\ba}_\kappa,$
\sn
\item[(b)]  $\mathbf a$ is a $T_{\ord}-(\theta_1,\theta_2)$-moral
  problem in $\mathbf m_0$ not solved in it.
\end{enumerate}
\end{enumerate}

[Why?  By \cite[Ch.VI,\S3]{Sh:c} there  is an ultra-filter $D$ on
$\lambda$ such that in $(\bbQ <)^\lambda/D$ there is a
 $(\theta_1,\theta_2)$-cut.  Define $\mathbf m$ by $\gB_{\mathbf m} =
  \cP(\lambda),D_{\mathbf m} = D$, now check. E.g.
   as $ \kappa = \lambda ^+ $ , easily the Boolean algebra
    $ \mathfrak{B} _ \mathbf{m} $ satisfies the 
    $ \kappa $-c.c.; 
    alternatively let $ \mathfrak{B} _ \mathbf{n} $ be generated
    by $ \{\mathbf{a} _{s,t}: s \in I_1, t \in I_2 \} $
    freely; and let $ D_ \mathbf{n} $ be the ultra-filter 
    on $ \mathfrak{B} _ \mathbf{n} $ generated by
    $\{  \mathbf{a} _{s,t}: s \in I_1, t \in I_2 \} $.
Now check.]

Let $\langle W_\alpha:\alpha < 2^\lambda\rangle$ be a partition of
$2^\lambda$ to sets each of cardinality $2^\lambda$ such that
$W_\alpha \cap \alpha = \emptyset$.

\begin{enumerate}
\item[$(*)_2$]  we can choose $\mathbf m_\alpha$ and $\langle \bar{\mathbf
  a}_\gamma:\gamma \in W_\alpha\rangle$ by induction on $\alpha \le
2^\lambda$ such that:
\sn
\begin{enumerate}
\item[$(a)$]  $\mathbf m_\alpha \in K^{\ba}_\kappa$ has cardinality $\le
  2^\lambda,$
\sn
\item[$(b)$]  $\langle \mathbf m_\beta:\beta \le \alpha\rangle \in
  S^{\ba}_\kappa,$
\sn
\item[$(c)$] $\mathbf m_0$ is as in $(*)_1,$
\sn
\item[$(d)$]  $\langle \bar{\mathbf a}_\gamma:\gamma \in
  W_\alpha\rangle$ 
be such that $\bar{\mathbf a}_\gamma$ is a
  $T_{\ord}-(\kappa_{\gamma,1},\kappa_{\gamma,2})$ problem in 
  $\mathbf m_\alpha$ and $\kappa_{\gamma,1},\kappa_{\gamma,2}$ are
regular $\le \lambda$ and $\{\theta_1,\theta_2\} \nsubseteq
\{\kappa_{\gamma,1},\kappa_{\gamma,2}\}$ and any such 
$\bar{\mathbf a}$ appears in the sequence,  
\sn
\item[$(e)$]  if $\alpha = \gamma +1$ then 
necessarily 
$\gamma \in W_\beta$ for
  some $\beta \le \alpha$ and in $\mathbf m_\alpha$ there is a solution
  for $\bar{\mathbf a}_\gamma,$
\sn
\item[$(f)$]  in $\mathbf m_\alpha$ there is no solution to $\bar{\mathbf
  a}^*$.
\end{enumerate}
\end{enumerate}

[Why can we carry the induction?

Now for $\alpha = 0$ use $(*)_1$, for $\alpha$ limit
use \ref{b13}
and for $\alpha$ successor use \ref{b20}.]

\begin{enumerate}
\item[$(*)_3$]  letting $\mathbf m = \mathbf m_{2^\lambda}$ we have
  $\gB_{\mathbf m} = \cup\{\gB_{\mathbf m_\alpha}:\alpha < 2^\lambda\}$ and
  $D_{\mathbf m} = \cup\{D_{\mathbf m_\alpha}:\alpha < 2^\lambda\}$.
\end{enumerate}

[Why?  Because $\langle \mathbf m_\alpha:\alpha \le 2^\lambda\rangle \in
S^{\ba}_\kappa$ and $\cf(2^\lambda) \ge \kappa$.]

\begin{enumerate}
\item[$(*)_4$]  there is a regular excellent filter $D_0$ on $\lambda$
  and homomorphism $\mathbf j$ from $\cP(\lambda)$ onto $\gB_{\mathbf m}$.
\end{enumerate}

[Why?  See \cite{Sh:997}.]

\begin{enumerate}
\item[$(*)_5$]  let $D = \mathbf j^{-1}(D_{\mathbf m})$.
\end{enumerate}

So $D$ is a filter on $\lambda$, and by \ref{g15} for $\theta =
\aleph_0$ (or Malliaris-Shelah \cite{Sh:997}) we are done.
%
\end{PROOF}

\begin{conclusion}
\label{b27}
If $\lambda \ge \aleph_2$ 
 \then \, 
the results of Malliaris-Shelah
\cite{Sh:998} cannot be generalized to reduced powers (atomic 
types, of course), that is (clause (A) is in contrast to
\cite[Th.10.25(b)$\Rightarrow $ (d)]{Sh:998}; 
clause (B) is in contrast to 
\cite[Th.10.1]{Sh:998}, and clause (C) is in contrast to
\cite[Th.3.1]{Sh:998})
\begin{enumerate}
\item[(A)] If $ \lambda \ge { \aleph_1 } $ then for some regular
 filter $ D $ on $ \lambda $ we have: in ultra-powers of 
 infinite linear orders
 we have a pre-cut with small cofinalities, but no symmetric pre-cut,
 that is:
 \begin{enumerate} 
 \item[(a)]  in the ultra-power $ (\mathbb{Q} , <  )^\lambda /D $ 
 there is a $ ({\aleph_1}, {\aleph_0}  ) $-pre-cut,
 
 \item[(b)] in this ultra-power, there is no symmetric pre-cut 
 of cofinality $ \sigma  $ for $ \sigma \le \lambda.$
 \end{enumerate} 
 
 \item[(B)] treetops: we can add above  that   in 
 $ ({}^{ \omega > } \omega , \triangleleft ) ^\lambda /D$
 every increasing sequence of length $ \le \lambda $ has an upper bound;
\item[(C)] if $ \lambda \le { \aleph_2 } $ 
then we can add in part (A),  
there are two pre-cuts
with the same small left cofinality but a different small right
cofinalities, e.g. $ {\aleph_1} $ from the left and $ { \aleph_2 }, \aleph_0$ from the right.
\end{enumerate}
\end{conclusion}

\begin{PROOF}{\ref{b27}}
For clause (A)  
we apply  clause (a) of \ref{b23}  choosing the pair 
$(\theta_1,\theta_2)$ as $(\aleph_1,\aleph_0)$.

For clause (B)   apply clause (b) of \ref{b23}.

For clause (C) we repeat the proof of \ref{b23} but starting
(with $ \kappa = \lambda ^+$ as there) 
and  choose as there $
\mathbf{m} _0 \in K_\kappa $ of cardinality 
$ \le 2^ \lambda $ such that some 
$ ( {\aleph_1}, {\aleph_0} )$-moral problem
and $ ( {\aleph_1},{ \aleph_2 } ) $-moral problem 
in $ \mathbf{m} _0 $ are not solve. Then continue as there.
\end{PROOF}



\begin{observation}
\label{b31}
If $\mathbf m \in K^{\ba}_\kappa$ \then \, any
$T_{\ord}-(\aleph_0,\aleph_0)$-problem
$\bar{\mathbf a}$  in $ \mathbf{m} $ has a solution.
\end{observation}

\begin{PROOF}{\ref{b31}}
Let $\mathbf b_{(1,n)} = \mathbf b_{(2,n)} = \mathbf b_n := \cap \{\mathbf
a_{(1,\ell),(2,k)}:\ell,k \le n\}$, clearly $s \in
I(\aleph_0,\aleph_0) \Rightarrow \mathbf b_s \in D$ and $(s,t) \in
I(1,\aleph_0) \times I(2,\aleph_0) \Rightarrow \mathbf b_s \cap \mathbf
b_t \le \mathbf a_{s,t}$.
\end{PROOF}

\begin{claim}
\label{b34}
In $M^{\gB}_*/D$, any increasing sequence of length $< \kappa^+$ has
an upper bound  \when \, (A) or (B) holds, where:

\begin{enumerate}
\item[$(A)$]  $(a) \quad M_* = ({}^{\omega >}\mu,\trianglelefteq),$
\sn
\item[${{}}$]  $(b) \quad \gB$ is a complete Boolean algebra which is
  $(< \theta)$-distributive,
\sn
\item[${{}}$]  $(c) \quad D$ is a $(\mu,\theta)$-regular,
  $\theta$-complete filter on $\gB,$
\sn
\item[${{}}$]  $(d) \quad (\bbQ,<)^{\gB}/D$ has no
  $(\sigma,\sigma)$-pre-cut for any regular $\sigma \le \kappa,$
  \sn
  \item[${{}}$]  $ ( e) \quad \mathbf{m} = (\mathfrak{B} , D ). $
\sn
\item[$(B)$]  $(a)-(c) \quad$ as above,
\sn
\item[${{}}$]  $(d) \quad$ every $T_{\tr}-(\sigma,\sigma)$-moral
  problem in $\mathbf m$ has a $T_{\tr}-(\sigma,\sigma)$-moral,

\hskip25pt   solution in $\mathbf m$ where:
\sn
\begin{enumerate}
\item[${{}}$]  $(\alpha) \quad \bar{\mathbf a}$ is a $T_{\tr}$-moral
  problem when:
\sn
\item[${{}}$]  \hskip10pt \quad $\bullet \quad \bar{\mathbf a} = \langle
  \mathbf a_{\alpha,\beta}:\alpha < \beta < \sigma\rangle,$
\sn
\item[${{}}$]  \hskip10pt \quad $\bullet \quad \mathbf a_{\alpha,\beta}
  \in D,$
\sn
\item[${{}}$]  \hskip10pt \quad $\bullet \quad$ if $u \subseteq
  \sigma$ is finite and $\mathbf c \in \gB^+$ then for some $\bar\eta = \langle
\eta_\alpha:\alpha \in u\rangle$ 

\hskip30pt  we have $\eta_\alpha \in
{}^{|u|>}|u|$ for $\alpha \in u$ and
$\mathbf c \le \mathbf a_{\alpha,\beta} \Rightarrow \eta_\alpha 
\trianglelefteq \eta_\beta$

\hskip30pt and $\mathbf c \cap \mathbf a_{\alpha,\beta} = 0_{\gB} \Rightarrow
\neg(\eta_\alpha \trianglelefteq \eta_\beta)$ for $\alpha < \beta$
from $u.$
\sn
\item[${{}}$]  $(\beta) \quad \bar{\mathbf b} = \langle \mathbf
  b_\alpha:\alpha < \sigma\rangle$ is a $T_{\tr}-
  (\sigma \sigma)$-solution of
  $\bar{\mathbf a}$ \when \, $\mathbf b_\alpha \in D$ and

\hskip25pt  $\mathbf b_\alpha \cap \mathbf b_\beta \le 
\mathbf a_{\alpha,\beta}$ for $\alpha < \beta < \sigma$.
\end{enumerate}
\end{enumerate}
\end{claim}

\begin{PROOF}{\ref{b34}}
If clause (A), as in \cite[Ch.VI,2.7]{Sh:a} or \cite{Sh:998}.

If clause (B), as above.
\end{PROOF}

\newpage

\bibliographystyle{amsalpha}
\bibliography{shlhetal}

\end{document}